\documentclass[11pt]{amsart}
\usepackage{amsmath,amssymb,amsfonts,url}
\numberwithin{equation}{section}

\newtheorem{thm}{Theorem}[section]
\newtheorem{lem}{Lemma}[section]
\newtheorem{rem}{Remark}[section]
\newtheorem{prop}{Proposition}[section]

\begin{document}
\title[Uniform estimate for Liouville system]{Profile of bubbling solutions to a Liouville system}

\author{Chang-shou Lin}
\address{Department of Mathematics\\
        Taida Institute of Mathematical Sciences\\
        National Taiwan University\\
         Taipei 106, Taiwan } \email{cslin@math.ntu.edu.tw}

\author{Lei Zhang}
\address{Department of Mathematics\\
        University of Alabama at Birmingham\\
        Birmingham, Alabama 35205}
\email{leizhang@math.uab.edu}
\thanks{Zhang is supported in part by NSF Grant 0600275 (0810902)}

\date{\today}

\begin{abstract}
In several fields of Physics, Chemistry and Ecology, some models
are described by Liouville systems. In this article we first prove
a uniqueness result for a Liouville system in $\mathbb R^2$. Then
we establish an uniform estimate for bubbling solutions of a
locally defined Liouville system near an isolated blowup point.
The uniqueness result, as well as the local uniform estimates are
crucial ingredients for obtaining a priori estimate, degree
counting formulas and existence results for Liouville systems
defined on Riemann surfaces.

\end{abstract}


\maketitle

\section{Introduction}

 In this article we are concerned with the following generalized Liouville
 system:
 \begin{equation}\label{1028e1}
 \Delta
u_i+\sum_{j=1}^na_{ij}h_je^{u_j}=0,\quad i\in I\equiv \{1,..,n\},
\quad \Omega\subset \mathbb R^2,
\end{equation}
where $\Omega$ is a subset of $\mathbb R^2$, $h_1,..,h_n$ are
positive smooth functions, $A=(a_{ij})_{n\times n}$ is an
invertible, symmetric and non-negative matrix .

(\ref{1028e1}) is an extension of the well known classical Liouville
equation
$$\Delta u+V e^u=0,\quad \Omega\subset \mathbb R^2,$$
which finds applications in many fields in Physics and
Mathematics. For example the Liouville equation is related to
finding a metric whose Gauss curvature is a prescribed function
\cite{chang}. In Physics, the Liouville equation represents the
electric potential induced by the charge carriers in electrolytes
theory \cite{rubinstein} and the Newtonian potential of a cluster
of self-gravitation mass distribution
\cite{aly,biler,wolansky1,wolansky2}. Moreover, it is closely
related to the abelian model in the Chern-Simons theories
\cite{jw1,jw2,jlw}.

 The Liouville systems are natural
extensions of the Liouville equation and they also have
applications in different fields of Physics, Chemistry and
Ecology. Indeed, various Liouville systems are used to describe
models in the theory of chemotaxis \cite{childress,keller},
in the physics of charged particle beams
\cite{bennet,debye,kiessling2} and in the theory of
semi-conductors \cite{mock}. For applications of Liouville
systems, see \cite{chanillo2,CSW} and the references therein. Here
we also note that another important extension of the Liouville
equation is the Toda system, which is closely related to the
non-abelian Chern-Simons theory \cite{dunne,yang}.

Chanillo and Kiessling \cite{chanillo2} first studied the type of
Liouville systems described by (\ref{1028e1}) with constant
coefficients in $\mathbb R^2$ and they proved that under certain
assumptions on $A$, all the entire solutions ($\Omega=\mathbb
R^2$) are symmetric with respect to some point. Their result was
improved by Chipot-Shafrir-Wolansky \cite{CSW}, who proved among
other things the following symmetry result:

\emph{ {\bf Theorem A} (Chipot-Shafrir-Wolansky) Let
$A=(a_{ij})_{n\times n}$ be a
\begin{equation}\label{matrix}
\mbox{invertible, symmetric, non-negative and irreducible matrix,}
\end{equation}
 $u=\{u_1,..,u_n\}$ be an entire solution of
\begin{equation}\label{chipotg}
\left\{\begin{array}{ll}
\Delta u_i+\sum_{j=1}^na_{ij}e^{u_j}=0,\quad \mathbb R^2,\\
\\
\int_{\mathbb R^2}e^{u_i}<\infty, \quad i\in I\equiv \{1,..,n\}.
\end{array}
\right.
\end{equation}
Then there exists $p\in \mathbb R^2$ such that all $u_1,..,u_n$
are radially symmetric and decreasing about $p$.}

\bigskip


 Recall that a matrix $A$ is called non-negative if $a_{ij}\ge
0$ ($i,j\in I$), irreducible if there is no partition of $I=I_1\cup
I_2, (I_1\cap I_2=\emptyset)$ such that $a_{ij}=0,\forall i\in
I_1,\forall j\in I_2$.

It turns out that the following quadratic polynomial is important
to the study of (\ref{chipotg}):
\begin{equation}\label{pisys}
\Lambda_J(\sigma)=4\sum_{i\in J}\sigma_i-\sum_{i,j\in
J}a_{ij}\sigma_i\sigma_j,\quad J\subset I\equiv \{1,..,n\},
\end{equation}
where $\sigma_i=\frac 1{2\pi}\int_{\mathbb R^2}e^{u_i}$,
$\sigma=\{\sigma_1,..,\sigma_n\}$.

 It was first proved by Chanillo-Kiessling \cite{chanillo2}
that entire solutions of (\ref{chipotg}) must satisfy a
Rellich-Pohozaev identity:
\begin{equation}\label{chanillopi}
\Lambda_I(\sigma)=4\sum_{i\in I}\sigma_i-\sum_{i,j\in I}
a_{ij}\sigma_i\sigma_j=0.
\end{equation}
Later Chipot-Shafrir-Wolansky \cite{CSW} proved the necessary and
sufficient condition for the existence of entire solutions to
(\ref{chipotg}):

\bigskip

\emph{{\bf Theorem B} (Chipot-Shafrir-Wolansky) Let $A$ satisfy
(\ref{matrix}). Then $\sigma=\{\sigma_1,..,\sigma_n\}$
satisfies
\begin{equation}\label{chipot3}
\Lambda_I(\sigma)=0,\,\, \mbox{and}\,\, \Lambda_J(\sigma)>0,\quad
\forall \emptyset\varsubsetneqq J\varsubsetneqq I,
\end{equation}
if and only if there exists a solution $\{u_1,..,u_n\}$ of
(\ref{chipotg}) such that $\frac 1{2\pi}\int_{\mathbb
R^2}e^{u_i}=\sigma_i$, $i\in I$.}

\bigskip

From now on we use $\Pi$ to represent the hyper-surface that satisfies (\ref{chipot3}).
It is immediate to observe that for each
$\sigma=\{\sigma_1,..,\sigma_n\}$ on $\Pi$, there
is more than one solution corresponding to $\sigma$. Indeed, let
$\{u_1,..,u_n\}$ be such a solution, then $\{v_1,..,v_n\}$ defined
by
$$v_i(y)=u_i(x_0+\delta y)+2\log \delta, \quad \forall x_0\in
\mathbb R^2, \quad \forall \delta>0 \quad i\in I$$ clearly solves
(\ref{chipotg}) and satisfies $\int_{\mathbb
R^2}e^{v_i}=\int_{\mathbb R^2}e^{u_i}$ ($i\in I$). A natural
question is: are all the solutions corresponding to $\sigma$
obtained from $\{u_1,..,u_n\}$ by translations and scalings? Our
first result in this paper is to give an affirmative answer to
this question:

\begin{thm}\label{thm6}
Let $A$ satisfy (\ref{matrix}),
$u=(u_1,..,u_n)$ and $v=(v_1,..,v_n)$ be two radial solutions of
(\ref{chipotg}) such that $\int_{\mathbb R^2}e^{u_i}=\int_{\mathbb
R^2}e^{v_i}$, $i\in I$, then there exists $\delta>0$ such that
$v_i(y)=u_i(\delta y)+2\log \delta$, $i\in I$.
\end{thm}

As is well known, for various equations it is important to have a
classification of all the global solutions. The classification
theorems of Caffarelli-Gidas-Spruck \cite{CGS},
Chen-Li\cite{chenli}, Jost-Wang \cite{jw1} and Lin \cite{lincla}
play a centrol role in the blowup analysis for prescribing scalar
curvature equations, prescribing Gauss curvature equations, Toda
systems and prescribing $Q-$curvature equations respectively.  The
existence result of Chipot-Shafrir-Wolansky (Theorem B) and the
uniqueness result (Theorem \ref{thm6}) can be combined to serve as
a classification theorem for the study of the blowup phenomena of
Liouville systems.

In \cite{CSW} Chipot-Shafrir-Wolansky also studied the Dirichlet
problem for the Liouville system (\ref{1028e1}) on bounded
domains. They considered the nonlinear functional $F$:
 $$F(u)=\frac 12\sum_{i,j\in I}\int_{\Omega}a^{ij}\nabla
 u_i\nabla u_j-\sum_{j\in I}\rho_j\log (\int_{\Omega}h_je^{u_j}),\quad u\in H_0^1(\Omega)$$
 where $a^{ij}$ ($i,j\in I$) are the entries of $A^{-1}$,
 $\rho_i$ ($i\in I$) are constants, and $h_i (i\in I)$ are positive
 smooth functions. Suppose the matrix $A=(a_{ij})$ is positive definite, it was shown in \cite{CSW} that $F$ is bounded
 from below in $H_0^1(\Omega)$ if and only if $\Lambda_I(\rho)\ge
 0$ ($\rho=(\rho_1,..,\rho_n)$), and a minimizer of $F(u)$ exists if
 $\Lambda_I(\rho)>0$. Obviously the Euler-Lagrange equation for
 the functional $F$ is the following

\begin{equation}\label{dirichlet}
\left\{\begin{array}{ll} \Delta
u_i+\sum_{j=1}^na_{ij}\rho_j\frac{h_je^{u_j}}{\int_{\Omega}h_je^{u_j}}=0,&\quad
\Omega\subset \mathbb R^2,
\quad i\in I\\
\\
u_i=0&\quad \mbox{on}\quad\partial\Omega
\end{array} \right.
\end{equation}
so the existence problem for (\ref{dirichlet}) is solved if
$\Lambda_I(\rho)>0$.

It is also natural to consider Liouville systems on Riemann
surfaces. Let $(M,g)$ be a Riemann surface of volume equal to $1$,
then the following variational form
$$J_{\rho}(u)=\frac 12\sum_{i,j=1}^n\int_{M}a^{ij}\nabla_g
u_i\nabla_g u_j+\sum_{j=1}^n\int_M\rho_ju_j-\sum_{j=1}^n\rho_j\log
\int_{M}h_je^{u_j} $$ corresponds to the system
\begin{equation}\label{1218e1}
\Delta_gu_i+\sum_{j=1}^n\rho_ja_{ij}(\frac{h_je^{u_j}}{\int_Mh_je^{u_j}dV_g}-1)=0,\quad
M, \quad i\in I.
\end{equation}

(\ref{dirichlet}) and (\ref{1218e1}) are generalizations of the
Liouville equation defined locally or on Riemann surfaces,
respectively. For the single Liouville equation, various results on
a priori estimate, degree counting formula and the existence of
solutions have been obtained by Chen-Lin \cite{ChenLin1,chenlin2}.
To study (\ref{dirichlet}) and (\ref{1218e1}), it is important to
understand the asymptotic behavior of blowup solutions.

 In this
article, we consider the following local estimate crucial to the
study of (\ref{dirichlet}) and (\ref{1218e1}): Let
$u^k=\{u_1^k,..,u_n^k\}$ be a sequence of functions which satisfies
\begin{equation}\label{assume1}\left\{\begin{array}{ll}
 \Delta
u_i^k+\sum_{j=1}^na_{ij}h_j^ke^{u_j^k}=0,\quad B_1\subset \mathbb R^2, \quad i\in I \\ \\
\int_{B_1}h_i^ke^{u_i^k}\le C,\,\, i\in I,\quad k=1,2,..
\end{array}
\right.
\end{equation}
where $B$ is the unit ball with center $0$, $\{h_i^k\}_{i\in I}$ are
positive $C^1$ functions uniformly bounded away from $0$:
\begin{equation}\label{assume4}
c_1^{-1}\le h_i^k\le c_1, \,\, \max_{B_1}|\nabla h_i^k|\le c_1,\,\,
i\in I,\quad k=1,2,...
\end{equation}
Suppose $0$ is the only blow-up point for $u^k$ and each component of
$u^k$ has a finite oscillation on $\partial B_1$:
\begin{eqnarray}\label{assume2}
\max_{\Omega}u_i^k\le C(\Omega),\quad \forall\,
\Omega\subset\subset B_1\setminus \{0\},\,\, i\in I\quad
k=1,2..\\
 \label{assume3} |u_i^k(x)-u_i^k(y)|\le c_0,\quad
\forall x,y\in
\partial B_1, \quad i\in I.
\end{eqnarray}

Our main assumption on $u^k$ is that $u^k$ converges to
a Liouville system of $n$ equations after scaling: Let $u_1^k(x_1^k)=\max_{B_1}u_i^k$
($i\in I$), $\epsilon_k=e^{-\frac 12u_1^k(x_1^k)}$ and
\begin{equation}\label{vik}
v_i^k(y)=u_i^k(\epsilon_ky+x_1^k)-u_1^k(x_1^k),\quad
y\in\Omega_k,\,i\in I
\end{equation}
 where $\Omega_k:=\{y;\,e^{-\frac
12u_1^k(x_1^k)}\cdot+x_1^k\in B_1\}$. Then
\begin{equation}\label{localcon}
 v^k=(v_1^k,..,v_n^k)\mbox{ converges in } C^2_{loc}(\mathbb
 R^2)\mbox{ to } v=(v_1,..,v_n)
 \end{equation}
 which is a solution of the Liouville system
 $$\Delta v_i+\sum_{j=1}^na_{ij}h_je^{v_j}=0,\,\, \mathbb R^2,\qquad
 h_i=\lim_{k\to \infty}h_i^k(x_1^k),\,\,i\in I.
 $$
 Note that $v_1,..,v_n$ are all radial functions because by Theorem
 A they are all radially symmetric with respect to a common point
 and $0$ is the maximum of $v_1$.
Our major local uniform estimate is:
 \begin{thm}\label{thmlocal} Let $A$ satisfy (\ref{matrix}), $u^k=(u_1^k,..,u_n^k)$
 be a sequence of solutions to (\ref{assume1}) such that
 (\ref{assume1})-(\ref{localcon}) hold. Then
 \begin{enumerate}
 \item there exists a sequence of radial solutions
 $V^k=(V_1^k,..,V_n^k)$ of
$$\Delta
V_i^k+\sum_{j=1}^na_{ij}h_j^k(0)e^{V_j^k}=0,\,\, \mathbb
R^2,\,\,\int_{\mathbb R^2}e^{V_i^k}<\infty,\,\, i\in I
$$
such that along a subsequence
\begin{equation}\label{uniformes}
|u_i^k(x)-V_i^k(x-x_1^k)|\le C(A,c_0,c_1,\sigma),\,\, i\in I\,\, x\in B_1,
\end{equation}
 where $\sigma=(\sigma_1,..,\sigma_n)$, $\sigma_i=\frac
1{2\pi}\int_{\mathbb R^2}h_ie^{v_i}$, $V^k$ is uniquely
determined by
\begin{enumerate}\item $V_1^k(0)=u_1^k(x_1^k)$
 \item $\int_{\mathbb
 R^2}h_j^k(0)e^{V_j^k}=\int_{B_1}h_j^ke^{u_j^k}$, \quad
 $j=1,..,n-1$.
 \end{enumerate}
\item There exists $\delta>0$ such that
$$\sum_{i,j\in I}a_{ij}\int_{B_1}h_i^ke^{u_i^k}\int_{B_1}h_j^ke^{u_j^k}
=8\pi\sum_{i\in I}\int_{B_1}h_i^ke^{u_i^k}+O(e^{-\delta
u_1^k(x_1^k)}).
$$
\end{enumerate}
\end{thm}

First we note that since every entire solution of the Liouville system satisfies
(\ref{chanillopi}),
$\int_{\mathbb R^2}h_n^k(0)e^{V_n^k}$ is
uniquely determined by (b) and
$$\sum_{ij}a_{ij}\int_{\mathbb R^2}h_i^k(0)e^{V_i^k}\int_{\mathbb R^2}h_j^k(0)e^{V_j^k}=8\pi\sum_i
\int_{\mathbb R^2}h_i^k(0)e^{V_i^k}. $$ Second, it is tempted to think that (\ref{uniformes}) is equivalent
to $|v_i^k-v_i|\le C$ ($i\in I$) in $\Omega_k$. In fact, the function $v$ may not be $V^k$ scaled
according to
the maximum of $u^k$ and the difference between $v^k$ and $v$ may not be uniformly bounded in $\Omega_k$. This is a
special feature of Liouville systems which can be observed from the entire solutions of
(\ref{chipotg}) as follows: Every point on $\Pi$ corresponds to
an entire solution. Let $\sigma^k=(\sigma^k_1,..,\sigma^k_n)$ be a sequence of points on
$\Pi$ that tends to $\sigma=(\sigma_1,..,\sigma_n)$. Let $\{w^k=(w^k_1,..,w^k_n)\}$ be a sequence of solutions
corresponding to $\sigma_k$ which converges in $C^2_{loc}(\mathbb R^2)$ to $w=(w_1,..,w_n)$, a solution
corresponding to $\sigma$. By standard potential analysis (see \cite{CSW})
$$w_i^k(x)=-(\sum_ja_{ij}\sigma_j^k)\ln |x|+O(1),\quad |x|>1$$ and
$$w_i(x)=-(\sum_ja_{ij}\sigma_j)\ln |x|+O(1),\quad |x|>1,\quad i\in I. $$ From the above we see that
 even though $\sigma^k\to \sigma$, the difference between
$w^k$ and $w$ may not be finite at infinity. Therefore the choice of $V^k$ in the statement of Theorem \ref{thmlocal}
is necessary.

For Liouville equations without singular data, the type of estimate
in Theorem \ref{thmlocal} was first derived by Li \cite{licmp}.
Later Bartolucci-Chen-Lin-Tarantello \cite{bclt} and Jost-Lin-Wang
\cite{jlw} established the same type of estimates for Liouville
equations with singular data and Toda systems, respectively. The
results of Li\cite{licmp} and Bartolucci-Chen-Lin-Tarantello
\cite{bclt} have been improved by Chen-Lin\cite{ChenLin1} and Zhang
\cite{zhangcmp,zhangccm} to a sharper form.

The estimates in Theorem \ref{thmlocal} would be very important when
a sequence of solutions $\{u^k\}$ of (\ref{1218e1}) has more than
one blowup point. Suppose $u^k=(u^k_1,..,u^k_n)$ is a sequence of
solutions of (\ref{1218e1}) with $\rho_i>0$ ($i\in I$). Assume that
$p_1,p_2$ are two blowup points, and the assumptions of Theorem
\ref{thmlocal} hold in neighborhoods around $p_1$ and $p_2$. By
Theorem \ref{thmlocal} there exist two entire solutions obtained
from the scaling of $u^k$ at $p_1$ and $p_2$. The question is
whether these two entire solutions are equal. Indeed, the answer is
yes when $A$ is positive definite, which is a consequence of Theorem
\ref{thm6} and \ref{thmlocal} (see section five for a proof of this
fact). The conclusion here is crucial to proving a priori estimates
for (\ref{dirichlet}) and (\ref{1218e1}). In a forthcoming paper
\cite{linzhang2} we shall discuss the a priori estimates, degree
counting formulas and existence results for (\ref{dirichlet}) and
(\ref{1218e1}).

Our next result concerns the location of blowup points for a
sequence of blowup solutions. Let $\{u^k\}$ be a sequence of
solutions of (\ref{assume1}) that satisfies the assumptions in
Theorem \ref{thmlocal}. Let $\{\psi_i^k\}_{i\in I}$ be the
harmonic functions
 defined by the oscillations of $u_i^k$ on $\partial B_1$:
$$\left\{\begin{array}{ll}
 \Delta \psi_i^k=0,\quad B_1,\\
 \\
 \psi_i^k=u_i^k-\frac 1{2\pi}\int_{\partial B_1}u_i^kdS, \quad
 \mbox{on}\quad \partial B_1,
 \end{array}
 \right.
 $$
 By the mean value property of harmonic functions we have
 $\psi_i^k(0)=0$. Also, since $\{u_i^k\}_{i\in I}$ have bounded oscillation on
 $\partial B_1$, all the derivatives of $\{\psi_i^k\}_{i\in I} $ on $B_{1/2}$
 are uniformly bounded.

\begin{thm}\label{location}
 Let $h_i$, $\psi_i$ ($i\in I$) be limits of $h_i^k$ and
 $\psi_i^k$, respectively, then under the same assumptions in
 Theorem \ref{thmlocal}
 $$\sum_{i\in I}(\frac{\nabla h_i(0)}{h_i(0)}+\nabla \psi_i(0))\sigma_i=0. $$
\end{thm}

Theorem \ref{location} can be used to determine the locations of
blowup points for (\ref{1218e1}) in the following typical situation. Let $\{u^k\}$ be a
sequence of blowup solutions to (\ref{1218e1}) with $\rho_i>0$
($i\in I$), $A$ satisfy (\ref{matrix}). In addition we assume $A$ to be positive definite
for simplicity. We can certainly assume $\int_Mh_i^ke^{u_i^k}dV_g=1$
($i\in I$) because for any solution $u=\{u_1,..,u_n\}$ to
(\ref{1218e1}), adding a constant vector $\{C_1,..,C_n\}$ to $u$
gives another solution. Suppose $p_1,..p_m$ are disjoint blowup
points of $u^k$ such that around each $p_t$ ($t=1,..,m$), $u^k$
converges in $C^2_{loc}(\mathbb R^2)$ to a Liouville system of $n$
equations after scaling. Let $G$ be the Green's function with
respect to $-\Delta_g$ on $M$:
$$
-\Delta_gG(x,p)=\delta_p-1, \quad \int_MG(x,p)dV_g(x)=0. $$
Corresponding to $G$ we define
$$G^*(x,p)=G(x,p)+\frac 1{2\pi}\chi(r) \log r$$
where $r=d_g(x,p)$, $\chi$ is a cut-off function supported in a
small neighborhood of $p$. Using $G^*$, the blowup points
$p_1,..,p_m$ are related by the following equation:
\begin{equation}\label{loca}
\sum_{i\in I}\bigg
(\frac{\nabla_gh_i(p_s)}{h_i(p_s)}+\frac{1}m(\sum_{j\in
I}\rho_ja_{ij})\sum_{t=1}^m\nabla_1G^*(p_s,p_t)\bigg )=0,\quad
s=1,..,m
\end{equation}
where $\nabla_1G^*$ means the covariant differentiation with
respect to the first component.

Even though the results in this paper (Theorems
\ref{thm6},\ref{thmlocal},\ref{location}) have their counterparts
for the Liouville equation, there are some essential differences
between the Liouville equation and the Liouville system that make
the analysis for the latter harder.
 First, the uniqueness theorem (Theorem \ref{thm6}) for
the system is generally harder to prove than one single equation,
because of the lack of the Sturm-Liouville comparison theory for the
linearized system. New ideas are needed to handle this difficulty.
In this article, we mainly use the method of continuation to prove
Theorem \ref{thm6}. Second, for the Liouville equation on $\mathbb
R^2$
$$\Delta u+e^u=0,\quad \mathbb R^2, \quad\int_{\mathbb
R^2}e^u<\infty. $$ All the solutions satisfy $\int_{\mathbb
R^2}e^u=8\pi$. However, for the Liouville system (\ref{chipotg}),
let $\sigma=(\sigma_1,..,\sigma_n)$ be the integration of the
entire solutions, which is on $\Pi$ (see (\ref{chipot3})). From Theorem B we
see that under some conditions we have a continuum of solutions,
as every point on $\Pi$ corresponds to a family of
solutions. This difference on the structure of entire solutions
exists not only between the Liouville equation and the Liouville
system, but also between the Liouville system and Toda systems
\cite{jlw}. Finally, for the Liouville equation, the Pohozaev
identity is a very useful tool, which gives a balancing condition
between the interior integration and the boundary integration.
However, for the Liouville system, the information from the
Pohozaev identity is limited, as we have more than one equation.
In this article, we use the uniqueness Theorem (Theorem
\ref{thm6}) to remedy what the Pohozaev identity can not provide.

The organization of the paper is as follows: In section two we prove Theorem
\ref{thm6} for two equations. We feel that the case of two equations is more explicit and represents
most of the difficulties of the system. Then in section three we prove the general case
of Theorem \ref{thm6} by mainly stating the difference with the proof in section two. In section four
we prove Theorem \ref{thmlocal} and in section five we prove Theorem \ref{location} as well as
(\ref{loca}). Finally in the appendix we list
a few Pohozaev identities to be used in different contexts.

\medskip

{\bf Acknowledgement}
 Part of the paper was finished when the second author was visiting Taida
 Institute of Mathematics during December 2007-January 2008 and in December 2008. He is
 very grateful to the Taida Institute of Mathematical Sciences for the warm
 hospitality.

\section{Proof of Theorem \ref{thm6} for two equations}

In this section we prove Theorem \ref{thm6} for two equations. So
the system is
\begin{equation}\label{july26e1}
\left\{\begin{array}{ll} \Delta u_1+a_{11}e^{u_1}+a_{12}e^{u_2}=0, \\
\Delta u_2+a_{12}e^{u_1}+a_{22}e^{u_2}=0,\quad \mathbb R^2,\\
\int_{\mathbb R^2}e^{u_1}<\infty, \quad \int_{\mathbb
R^2}e^{u_2}<\infty
\end{array}
\right.
\end{equation}
where the assumption on $A$ now becomes $a_{ii}\ge 0,\, i=1,2$, $a_{12}>0$ and
$a_{12}^2\neq a_{11}a_{22}$. Let
$$\sigma_i=\frac 1{2\pi}\int_{\mathbb R^2}e^{u_i}
\quad \mbox{and} \quad m_i=\sum_ja_{ij}\sigma_j, \quad i\in I=\{1,2\}. $$ By
standard potential analysis (see, for example \cite{CSW}) we have
\begin{equation}\label{mig2}
m_i>2\quad i\in I=\{1,2\}
\end{equation}
 and
\begin{equation}\label{uinfinity}
u_i(x)=-m_i\ln |x|+O(1),\quad |x|>1\quad i\in I.
\end{equation}

Let $u=\{u_1,u_2\}$ be a radial solution of (\ref{july26e1}) and we
consider the linearized equation of (\ref{july26e1}) at $u$:
\begin{equation}
\label{jun26e1} (r \phi_i'(r))'+\sum_ja_{ij}
e^{u_j}\phi_j(r)r=0,\quad 0<r<\infty,\quad i\in I.
\end{equation}

\begin{lem}\label{loginfinity}
Let $\phi=(\phi_1,\phi_2)$ be a solution of (\ref{jun26e1}), then
$\phi_i(r)=O(\ln r)$ at infinity for $i\in I$.
\end{lem}

\noindent{\bf Proof of Lemma \ref{loginfinity}:} Let
$\psi(t)=(\psi_1(t), \psi_2(t))$ be defined as
$$\psi_i(t)=\phi_i(e^t),\quad i\in I. $$
Then $\psi$ satisfies
$$\psi_i''(t)+\sum_ja_{ij}e^{u_j(e^t)+2t}\psi_j(t)=0,\quad
-\infty<t<\infty,\quad i\in I. $$ Let
$\psi_3=\psi_1'$, $\psi_4=\psi_2'$
and ${\mathbf F}=(\psi_1,..,\psi_4)^T$, then $\mathbf{F}$
satisfies
$$\mathbf{F}'=\mathbf{MF}$$

where $\displaystyle{\mathbf{M}=\left(\begin{array}{cc}
\mathbf{0}&\mathbf{I}\\
\mathbf{B}&\mathbf{0}
\end{array}
\right) }$. $\mathbf{B}$ is a $2\times 2$ matrix with
$\mathbf{B}_{ij}=-a_{ij}e^{u_j(e^t)+2t}$. For $t>1$, the solution
for $\mathbf{F}$ is
\begin{equation}\label{jan22e1}
\mathbf{F}(t)=\lim_{N\to \infty}e^{\epsilon \mathbf{M}(t_N)}...e^{\epsilon \mathbf{M}(t_0)}\mathbf{F}(0).
\end{equation}
where $t_0,...,t_N$ satisfy $t_j=j*\epsilon$, $j=0,..,N$, $\epsilon=t/N$.
 Since $u_i(e^t)+2t \sim (-m_i+2)t$ when $t$ is large and
$m_i>2$(see (\ref{mig2})\, ), we have $\|\mathbf{B}\|\sim
e^{-\delta t}$ for some $\delta>0$ and $t$ large. With this property we further have
\begin{equation}\label{jan22e2}
\|\mathbf{M}\|^k\le Ce^{-k\delta_1t},\quad k=2,3,... \quad t>0
\end{equation}
for some $\delta_1>0$.
Using (\ref{jan22e2}) in (\ref{jan22e1}) we have
$$\|\mathbf{F}(t)\|=O(t),\quad t>1. $$
Lemma \ref{loginfinity} is established. $\Box$

\begin{lem} \label{SetS}
Let $\phi=\{\phi_1,\phi_2\}$ be a bounded solution of
(\ref{jun26e1}), then $\phi=C(ru_1'+2,ru_2'+2)$ for some
constant $C$.
\end{lem}

\noindent{\bf Proof of Lemma \ref{SetS}:} Let
$$\phi^0=(ru_1'+2, ru_2'+2), $$
it is easy to verify that $\phi^0$ solves
(\ref{jun26e1}) and $\phi^0$ is bounded. We prove Lemma
\ref{SetS} by contradiction. Suppose $\bar \phi=(\bar \phi_1,\bar \phi_2)$ is
another bounded solution of (\ref{jun26e1}) and is not a multiple
of $\phi^0$, then $\phi^0$ and $\bar \phi$ form a basis for
all the solutions of (\ref{jun26e1}). Since $\bar \phi_1(0)$ and $\bar \phi_2(0)$ can not both be $2$,
without loss of generality
we assume $\bar \phi_1(0)=0$ and $\bar \phi_2(0)=1$. We use $E$ to
denote the set of all solutions. Since every solution is a linear
combination of $\phi^0$ and $\bar \phi$,
all the solutions are bounded. Let
\begin{eqnarray*}
 S=\{\alpha |\,\, \owns (\phi_1,\phi_2)\in E, \, \phi_1(0)=2,
\phi_2(0)=\alpha\le 2, \, \mbox{such that }\\
\int_0^re^{u_i}\phi_i(s)sds>0 \,\,
\mbox{for all } r>0\, \,\, i\in I
\}.
\end{eqnarray*}

We note that if $\phi_2(0)=2$, then
$\phi(r)=(ru_1'(r)+2,ru_2'(r)+2)$. It is easy to see that $2\in S$ because
$$
\int_0^re^{u_i}(su_i'+2)sds=r^2e^{u_i(r)}>0, \quad i\in I.
$$
 Next we see that $S$ is a bounded set. Because if
$\alpha<0$, let $\phi=\{\phi_1,\phi_2\}$ be the bounded solution such that
$\phi_1(0)=2,\phi_2(0)=\alpha$. Then $\int_0^re^{u_2(s)}\phi_2(s)sds<0$
for $r$ small enough. So $\alpha \not \in S$.

Set $\alpha_0=\inf_S\alpha$. Then we claim that $\alpha_0\in S$.
In fact, let $\{\alpha_k\in S\}$ tend to $\alpha_0$ from above as $k\to \infty$,
let $\phi^k=\{\phi^k_1,\phi^k_2\}$ correspond to $\alpha_k$. Since
$\alpha_k\in S$, $\int_0^rse^{u_i}\phi^k_i(s)ds>0$, for all $r$. Moreover, it is easy
to see that $\phi^k$ converge to a solution
$\phi=(\phi_1,\phi_2)$ in $E$ because $\phi^k$s are linear combinations
of $\phi^0$ and $\bar \phi$. It is also
immediate to observe from the convergence that
$$\int_0^re^{u_i}\phi_i(s)sds\ge 0, \quad
\mbox{for all }r>0,\quad i\in I.
$$
Thus,
$$
r\phi_i'(r)=-\sum_ja_{ij}\int_0^re^{u_j}\phi_j(s)sds\le
0, \quad i\in I.
$$
So both $\phi_1$ and $\phi_2$ are non-increasing functions. Since they are bounded functions, for each $i\in I$
there exist $r_l\to
\infty$ such that $r_l\phi_i'(r_l)\to 0$, which leads to
$$\sum_ja_{ij}\int_0^{\infty}e^{u_j}\phi_j(s)sds=0,\quad i\in I.
$$
Then we obtain the following from the invertibility of $A$:
\begin{equation}\label{phipsi}
\int_0^{\infty}e^{u_i}\phi_i(s)sds=0, \quad i\in I.
\end{equation}

Since $\phi_1$ and $\phi_2$ are non-increasing functions, (\ref{phipsi})
implies that
$$\lim_{r\to \infty}\phi_i(r)<0\quad i\in I.$$
Indeed, for example for $\phi_1$,
$\int_0^{\infty}e^{u_1}\phi_1(s)sds=0$ and the monotonicity or $\phi_1$
imply either $\lim_{r\to \infty}\phi_1(r)<0$ or $\phi_1\equiv 0$. Then
we see immediately that the latter case does not occur, as
$\phi_1(0)=2$. Similarly for $\phi_2$, the case that $\phi_2\equiv 0$ also
does not happen because $\phi_1\not \equiv 0$. Another immediate
observation is $\phi_2(0)>0$.

For the above, we have

$$
\int_{0}^{r} e^{u_i}\phi_{i}sds>0 \quad\mbox{if}\quad\phi_{i}(r)\geq
0,\quad\mbox{and}
$$

$$
\int_{0}^{r} e^{u_i}\phi_{i}sds>\int_{0}^{\infty}e^{u_i}\phi_{i}sds
=0, \quad\mbox{if}\quad\phi_{i}(r)< 0.
$$

\noindent Thus $\alpha_0\in S$.

Now we claim that for $\epsilon>0$ small enough,
$\alpha_0-\epsilon\in S$. Indeed, consider $\phi-\epsilon\bar \phi$,
obviously this is a solution to (\ref{jun26e1}) and satisfies
$\phi_1(0)-\epsilon \bar \phi_1(0)=2$, $\psi(0)-\epsilon \bar
\psi(0)=\alpha_0-\epsilon$. Since $\{\phi_1-\epsilon\bar
\phi_1,\phi_2-\epsilon\bar \phi_2\}$ is a bounded solution of
(\ref{jun26e1}) we have
$$\int_0^{\infty}e^{u_i}(\phi_i-\epsilon\bar \phi_i)sds=0, \quad
i\in I.$$ For $r$ large and $\epsilon$ small, since $\phi_1(r)$ and
$\phi_2(r)$ are smaller than a negative number for $r$ large, it is
easy to choose $\epsilon$ small enough so that

$$
\int_r^{\infty}e^{u_i}(\phi_i-\epsilon\bar \phi_i)sds<0,\quad i\in i
$$
for all large $r$ large. Consequently
\begin{equation}\label{mono}
\int_0^re^{u_i}(\phi_i-\epsilon\bar \phi_i)sds>0\quad i\in I
\end{equation}
for all large $r$. Then by possibly choosing $\epsilon>0$ smaller,
we can make (\ref{mono}) hold for all $r>0$. $\alpha_0-\epsilon\in
S$ is proved. This is a contradiction to the definition of
$\alpha_0$. Lemma \ref{SetS} is established. $\Box$

\medskip

Now we are in the position to complete the proof of Theorem
\ref{thm6} for two equations. We consider the following
initial-value problem:

\begin{equation}\label{july14e5}
\left\{\begin{array}{ll}
{u_i}^{''}+\frac{u_i'}{r}+\sum_ja_{ij}e^{u_j}=0,\quad i=1,2, \\
\\
u_1(0)=\alpha,\quad u_2(0)=0.
\end{array}
\right.
\end{equation}

\medskip

\noindent Case 1: $a_{ii}>0, i=1,2$

\medskip

Since $a_{ii}>0$, by Lemma 3.2 in section three, the solution pair
$u_i(r)$, exists for all $r>0$ and $i=1,2$, and satisfies

$$
\int_{0}^{\infty} e^{u_i(r)} rdr < +\infty \quad, i=1,2.
$$

Set
$$
\sigma_i(\alpha)=\int_{0}^{\infty} e^{u_i(r)}rdr \quad, i=1,2.
$$

Thus $\sigma(\alpha)=(\sigma_1,\sigma_2)$ is a function of $\alpha$
and lies in $\Pi$ (defined by (\ref{chipot3})), which is a curve: $\Lambda_I(\sigma)=0$ ($\sigma_1,\sigma_2>0$).
We want to prove that
$$
\sigma: \quad \mathbb R \quad\rightarrow\quad
\Pi
$$
is an 1-1 and onto map. Since both $\mathbb R$ and
$\Pi$ are connected, it suffices to prove
$\sigma$ is an open mapping. In the following, we want to show the
claim
\begin{equation}\label{july14e6}
\frac{\partial\sigma_1}{\partial\alpha}\neq 0 \quad\mbox{and}\quad
\frac{\partial\sigma_2}{\partial\alpha}\neq 0 \quad\mbox{for
all}\quad\alpha\in\mathbb R^2
\end{equation}
Then the openness of $\sigma$ follows immediately.

We prove this claim by contradiction. Suppose there exists $\alpha$
such that, say, $\partial_{\alpha}\sigma_1=0$. This implies
immediately that

\begin{equation}\label{jan23e1}
\int_0^{\infty}re^{u_1}\phi_1=0,
\end{equation}

\noindent where $\phi_1=\partial_{\alpha}u_1$. Correspondingly we
set $\phi_2=\partial_{\alpha}u_2$. Then $\{\phi_1,\phi_2\}$
satisfies the linearized system (\ref{jun26e1}). By Lemma
\ref{loginfinity} $\phi_i(r)=O(\ln r)$ at infinity. The Pohozaev
identity for (\ref{jun26e1}) is (see the appendix for the proof)
\begin{equation}\label{boundedatinfi}
\sum_i(r^2\phi_i(r)e^{u_i}-2\int_0^rse^{u_i}\phi_i(s)ds)
=-\sum_{ij}a^{ij}(r\phi_i'(r))(ru_j'(r)).
\end{equation}
The first term on the left hand side of (\ref{jan23e1}) tends to $0$
as $r\to \infty$. To deal with the terms on the right hand side,
first we use the equation for $\phi_i$ to get
$$-r\phi_i'(r)=\sum_la_{il}\int_0^rse^{u_l}\phi_l(s)ds. $$
The equation for $u_i$ gives $\lim_{r\to \infty} ru_i'(r)=-m_i$.
Putting the above information together we obtain the following from
(\ref{boundedatinfi}):
$$\sum_i(m_i-2)\int_0^{\infty}se^{u_i}\phi_i(s)ds=0. $$
By (\ref{jan23e1}), we have
$$
\int_0^{\infty}e^{u_i}\phi_i rdr=0,\quad i=1,2.
$$
Using (\ref{boundedatinfi}) for the equation for $\phi_i$ we have
$$-r\phi_i'(r)=\int_0^r\sum_ja_{ij}e^{u_j}\phi_jsds=-\int_r^{\infty}\sum_ja_{ij}e^{u_j}\phi_jsds
=O(r^{-\delta})$$ for some $\delta>0$. Therefore $\phi_i$ ($i\in I$)
is bounded at infinity.
 By Lemma \ref{SetS}, there is a constant $c$ such that $\phi_1=c(r{u_1}'+2),\phi_2=c(r{u_2}'+2)$.
But one sees immediately that this is impossible because
$\phi_1(0)=0$, $\phi_2(0)=1$. The claim is proved.

\medskip

Theorem \ref{thm6} for this case is implied by the claim. In fact, suppose $\{\bar
u_1,\bar u_2\}$ is another pair of radial solutions of the Liouville
system so that $\int_{\mathbb R^2}e^{\bar u_i}=\int_{\mathbb
R^2}e^{u_i}$ ($i=1,2$). By scaling, we may assume $u_2(0)=\bar
u_2(0)=0$. Since the mapping $\sigma: \mathbb R^{n-1}\to \Pi$ is one to one and onto,
we have $u_1(0)=\bar u_1(0)$. Consequently
$u_i\equiv \bar u_i$ ($i\in I$), hence Theorem \ref{thm6} is proved
for the case $a_{ii}>0,i=1,2$.

\medskip
\noindent Case 2: There exists $i$ such that $a_{ii}=0$.
\medskip

Set
$$
\Pi_1=\{\alpha|e^{u_j}\in L^1(\mathbb R^2),\quad j=1,2,\quad u=(u_1,u_2) \mbox{ is a solution of
(\ref{july14e5})}\}
$$
Similar to the previous step, the map
$\Pi_1\rightarrow \Pi$ is an open mapping.
Since $a_{11}=0$ or $a_{22}=0$, $\Pi$ is
non-compact and connected. Thus $\sigma$ is 1-1 and onto from each
component of $\Pi_1$ onto $\Pi$.

Now suppose $\Pi_1$ has two component, say $\Pi_1^1$ and $\Pi_1^2$. Choose any
$\sigma$ of $\Pi$. Then there exists $\alpha_1\in
\Pi^1_1$, and $\alpha_2\in \Pi_1^2$ such that $u^1=({u_1}^1,{u_2}^1)$ and
$u^2=({u_1}^2,{u_2}^2)$ are the corresponding solutions of (\ref{july14e5}) and
satisfy
$$
\int_{0}^{\infty} e^{{u_j}^1} rdr=\int_{0}^{\infty} e^{{u_j}^2}
rdr=\sigma_j\quad j=1,2.
$$
Clearly, $\exists R_0$ such that for $r\geq R_0$ and some
$\delta>0$,
$$
({u_j}^k)^{'}(r)r\leq -(2+2\delta)\quad j=1,2 \quad k=1,2.
$$

Now consider the perturbation of (\ref{july14e5}):

\begin{equation}\label{pertur1}
\left\{
\begin{array}{ll}
\Delta u_i+\sum_{j=1}^2(a_{ij}+\epsilon \delta_{ij}) e^{u_j}=0 \quad \mathbb R^2,\quad i=1,2\\
\\
u_1(0)=\alpha, \quad u_2(0)=0.
\end{array}
\right.
\end{equation}
Here we require $\epsilon\in (0,\delta_0)$ where $\delta_0$ is so small that the matrix $(a_{ij}+\epsilon\delta_{ij})_{n\times n}$
is non-singular for all $\epsilon\in (0,\delta_0)$.
Let $u^{k,\epsilon}=(u^{k,\epsilon}_1,u^{k,\epsilon}_2)$ be the solution of (\ref{pertur1}) with respect to the initial
condition $(\alpha_k,0)$ ($k=1,2$).
For $\delta_0$ small we have
$$
(u^{k,\epsilon}_j(r))'r\leq -(2+\delta)\quad\mbox{at}\quad
r=R_0,\quad 0\le \epsilon \le \delta_0.
$$
Then by the super-harmonicity of $u^{k,\epsilon}_j$ it is easy to show
$$
(u^{k,\epsilon}_j(r))'r\le -(2+\delta)\quad\mbox{for}\quad r\geq
R_0.
$$

\noindent Thus, $\exists C>0$ and $R_1\geq R_0$ such that
\begin{equation}\label{feb25e3}
e^{u^{k,\epsilon}_j(r)}\leq Cr^{-(2+\delta)}\quad\mbox{for}\quad
r\geq R_1
\end{equation}

Hence for $k=1,2$,
$$
{\sigma_j}^{\epsilon}(\alpha_k)=\int_0^\infty
e^{u^{k,\epsilon}_j(r)}rdr=\int_0^\infty
e^{u_j^k(r)}rdr+o(1)=\sigma_j+\circ(1), \quad j=1,2.
$$
where $o(1)\rightarrow 0$ as $\epsilon\rightarrow 0$.

Next we claim that
\begin{equation}\label{energyder}
\frac{\partial \sigma^{\epsilon}_j}{\partial \alpha}(\alpha_k)=\frac{\partial \sigma_j}{\partial \alpha}(\alpha_k)+\circ(1).
\end{equation}
Indeed,
\begin{equation}\label{feb25e1}
\frac{\partial \sigma^{\epsilon}_j}{\partial \alpha}(\alpha_k)=\int_0^{\infty}re^{u^{k,\epsilon}_j(r)}
\frac{\partial u^{k,\epsilon}_j}{\partial \alpha}(r)dr,\quad j=1,2,\quad k=1,2.
\end{equation}
$(\frac{\partial u^{k,\epsilon}_1}{\partial \alpha},\frac{\partial u^{k,\epsilon}_2}{\partial \alpha})$ satisfies the
following linearized equation:
$$-\Delta (\frac{\partial u^{k,\epsilon}_i}{\partial \alpha})=\sum_{j=1}^2(a_{ij}+\epsilon \delta_{ij})e^{u^{k,\epsilon}_j}
\frac{\partial u^{k,\epsilon}_j}{\partial \alpha},\quad i=1,2.
$$
Using the argument of Lemma \ref{loginfinity} we have
\begin{equation}\label{feb25e2}
|\frac{\partial u^{k,\epsilon}_i}{\partial \alpha}(r)|\le C\ln r,\quad r\ge 2, \quad i=1,2.
\end{equation}
where the constant $C$ is independent of $\epsilon\in (0,\delta_0)$. Moreover, for any fixed $R>0$,
$\frac{\partial u^{1,\epsilon}_i}{\partial \alpha}(r)$ converges uniformly to
$\frac{\partial u^1_i}{\partial \alpha}(r)$ over $0<r<R$ with respect to $\epsilon$.
Using the decay estimates (\ref{feb25e3}) and (\ref{feb25e2}) in (\ref{feb25e1}) we obtain (\ref{energyder})
by elementary analysis.

Since $\lim_{\epsilon\rightarrow
0}\frac{\partial{\sigma_j}^\epsilon}{\partial\alpha}(\alpha_1)
=\frac{\partial\sigma_j}{\partial\alpha}(\alpha_1)\neq 0$, there
exists $\alpha_1(\epsilon)=\alpha_1+o(1)$ such that
\begin{equation}\label{feb25e4}
{\sigma_1}^\epsilon(\alpha_1(\epsilon))={\sigma_1}^\epsilon(\alpha_2).
\end{equation}
Both $(\sigma_1^{\epsilon}(\alpha_1(\epsilon)), \sigma^{\epsilon}_2(\alpha_1(\epsilon)))$ and
$(\sigma_1^{\epsilon}(\alpha_2), \sigma^{\epsilon}_2(\alpha_2))$ satisfy
${\Lambda_I}^\epsilon(\sigma^\epsilon)=0$, which reads
$$\sum_{i,j=1}^2(a_{ij}+\epsilon \delta_{ij})\sigma^{\epsilon}_i\sigma^{\epsilon}_j=4\sum_{i=1}^2\sigma_i^{\epsilon}. $$
Using (\ref{feb25e4}) in the above we have
$$
{\sigma_2}^\epsilon(\alpha_1(\epsilon))={\sigma_2}^\epsilon(\alpha_2).
$$

Since $\alpha_1(\epsilon)\neq\alpha_2$, it yields a contradiction to
the uniqueness property that the system (\ref{pertur1}) satisfies. Hence the proof of
Theorem \ref{thm6} for two equations is complete.  $\Box$






\section{Proof of Theorem \ref{thm6} for the general case}

The proof for the general case of Theorem \ref{thm6} is similar to the case of two equations. We mainly focus
on the difference in this section.

First we point out that Lemma \ref{loginfinity} still holds for the general case with the same proof. The first major
result in this section is the following
\begin{lem}\label{jan20lem}
Let $\phi=(\phi_1,..,\phi_n)$ be a bounded solution of
\begin{equation}\label{jan20e1}
(r\phi_i'(r))'+\sum_{j=1}^na_{ij}e^{u_j}r\phi_j(r)=0,\quad 0<r<\infty, \quad i\in I=\{1,..,n\},
\end{equation}
then $\phi_i(r)=ru_i'(r)+2$, $i\in I=\{1,..,n\}$.
\end{lem}

\noindent{\bf Proof of Lemma \ref{jan20lem}:} Let $\phi^0=(ru_1'(r)+2,..,ru_n'(r)+2)$, then by direct computation one sees that
$\phi^0$ is a solution of (\ref{jan20e1}). Suppose there is another bounded solution $\phi^1=(\phi^1_1,..,\phi^1_n)$
different from $\phi^0$, without
loss of generality we assume $\phi_1(0)=0$, as one of $\phi_i^1(0)$ must be different from $2$. To derive
a contradiction we define
\begin{eqnarray*}
S=\{\alpha;\quad \owns \mbox{ a  bounded solution }\phi\mbox{ such that } \phi_1(0)=2,\\
\phi_i(0)=\alpha_i\le 3,\,\, i=2,..,n; \quad
\alpha=\min\{\alpha_2,..,\alpha_n\}\\
\int_0^re^{u_i(s)}\phi_i(s)sds>0,\,\, \forall r>0, \,\, i\in I\}.
\end{eqnarray*}

By direct computation $2\in S$, which corresponds to the solution $\phi^0$. Since $\phi^0_i$ ($i\in I$) is strictly decreasing, we can
choose $t$ small enough to make all components of
$\phi^0+t\phi^1$ strictly decreasing. By choosing $t$ or $-t$ we can make $2-\epsilon\in S$ for some $\epsilon>0$ sufficiently small.
Let $\bar \alpha$ be the infimum of $S$ and let $\alpha^k=\{\alpha^k_1,..,\alpha^k_n\}\in S$ be a sequence in $S$ that tends to $\bar \alpha$ from above.
Suppose $\phi^k=\{\phi_1^k,..,\phi_n^k\}$ is the solution corresponding to $\alpha^k$, then we claim that $\{\phi^k\}$ converges
to $\bar \phi=\{\bar \phi_1,..,\bar \phi_n\}$, which is also a bounded solution with strict monotone properties described in $S$.
Indeed, let $\psi^m=(\psi^m_1,..,\psi^m_n)$ be the solution to (\ref{jan20e1}) such that $\psi^m_j(0)=\delta_j^m$. By Lemma
\ref{loginfinity} $\psi^m_i(r)=O(\ln r)$ at infinity. $\phi^k$ can be written as
\begin{equation}\label{jan20e3}
\phi^k=\sum_{m=1}^n\alpha^k_m\psi^m.
\end{equation}
Since $\bar \alpha\le \alpha_i^k\le 3,$ ($i\in I$) for all $k$, along a subsequence, $\alpha^k$ converges to $\{\bar \alpha_1,
,,\bar \alpha_n\}$. As a consequence, $\phi^k$ converges to $\bar \phi=\sum_{m=1}^n\bar \alpha_m\psi^m$ uniformly over any compact
subsets of $\mathbb R^2$. The monotone property of $\phi^k$ implies that
$$\int_0^re^{u_i}\bar \phi_i(s)sds\ge 0, \quad i\in I, \quad \forall r>0. $$
On the other hand, since $\phi^k$ are all bounded functions, for each $\phi_i^k$ we find $r_l\to \infty$
such that $r_l(\phi_i^k)'(r_l)\to 0$. This leads to
$$\int_0^{\infty}\sum_ja_{ij}e^{u_j(s)}\phi_i^k(s)s ds=0 \quad i\in I.$$
Since $A$ is invertible we have
\begin{equation}\label{jan23e3}
0=\int_0^{\infty}e^{u_i}\phi_i^k(s)sds=\sum_{m=1}^n\alpha_m^k\int_0^{\infty}e^{u_i(s)}\psi_i^m(s)sds, \quad i\in I.
\end{equation}
Since $\int_0^{\infty}e^{u_i}\psi^m_i(s)sds$
is well defined, we let
$\alpha^k\to (\bar \alpha_1,..,\bar \alpha_n)$ to get
\begin{equation}\label{jan20e4}
\int_0^{\infty}e^{u_i(s)}\bar\phi_i(s)s=0,\quad i\in I.
\end{equation}
Using the argument for the case of two equations as well as the assumption that $A$ is irreducible we know
each $\bar \phi_i$ decreases into a negative constant at infinity and $\bar \phi_i(0)>0$.
As a consequence, $\int_0^re^{u_i(s)}\bar \phi_i(s)sds>0$ for each $r>0$ and $\bar \alpha>0$. Thus $\bar \alpha\in S$. Then as in
the case for two equations, $\{\bar \phi+t\phi^1\}$ for $t$ small enough also satisfies the strict monotone property
described in the definition of $S$. Therefore $\bar \alpha-\epsilon\in S$ for $\epsilon>0$ small enough. This is a contradiction
to the definition of $\bar \alpha$. Lemma \ref{jan20lem} is established. $\Box$

\medskip

Now we complete the proof of Theorem \ref{thm6} for $n$ equations.
Let $u=(u_1,..,u_n)$ satisfy
\begin{equation}\label{jan21e1}
\left\{\begin{array}{ll}
u_i''(r)+\frac{u_i'(r)}r+\sum_ja_{ij}e^{u_j}=0\quad 0<r<\infty,\quad i\in I\\
\\
\int_0^{\infty}re^{u_i(r)}dr<\infty, \\
\\
u_1(0)=\beta_1, \,\, ... \,\, ,u_{n-1}(0)=\beta_{n-1},\quad u_n(0)=0.
\end{array}
\right.
\end{equation}

The following Lemma is useful for the case $a_{ii}>0$.
\begin{lem}\label{domain}
Let $a_{ii}>0$ ($i\in I$), then for all $\beta=(\beta_1,..,\beta_{n-1})\in \mathbb R^{n-1}$, there
exists a solution $u=(u_1,..,u_n)$ to (\ref{jan21e1}).
\end{lem}

\noindent{\bf Proof of Lemma \ref{domain}:} By standard ODE
existence theory we see that for $\beta=(\beta_1,..,\beta_{n-1})\in
\mathbb R^{n-1}$, there exists a radial solution $u=(u_1,..,u_n)$ in
the neighborhood of $0$. Then by writing the system as a first order
ODE system we see the right hand side always satisfies the Lipschitz property,
therefore by Picard's theorem the solution exists for all $r>0$. We are left to
show that $\int_0^{\infty}e^{u_i(s)}sds<\infty$. Let $v_i(t)=
u_i(e^t)+2t$ ($i\in I$), then $v=(v_1,..,v_n)$ satisfies
$$v_i''(t)+\sum_ja_{ij}e^{v_j(t)}=0,\quad -\infty<t<\infty,\quad
i\in I.$$ From the equation for $u_i$ we have
$$ru_i'(r)=-\int_0^r\sum_ja_{ij}e^{u_j(s)}sds<0,\quad r>0,\quad i\in
I. $$ Consequently $v_i'(t)<2$ for $t\in \mathbb R$. Fix $t_0\in
\mathbb R$ we have, for $t>t_0$,
$$v_i'(t)=v_i'(t_0)-\int_{t_0}^t\sum_ja_{ij}e^{v_j(s)}ds,\quad i\in
I.$$ Since $a_{ii}>0$ and $a_{ij}\ge 0$, it is easy to see that
there exists $t>t_0$ such that $v_i'(t)<0$. Choose $t_1$ such that
$v_i'(t_1)=-\delta<0$ for some $\delta>0$, then we see from the equation for $v_i$ that
$$v_i(t)\le v_i(t_1)-\delta (t-t_1),\quad t>t_1 $$
which is equivalent to $u_i(r)<(-2-\delta)\ln r+C$ for $r>e^{t_1}$.
Therefore $\int_0^{\infty}e^{u_i(s)}sds<\infty$. Lemma
\ref{domain} is established. $\Box$

\bigskip

Recall that  $\sigma_i=\frac 1{2\pi}\int_{\mathbb
R^2}e^{u_i}=\int_0^{\infty}e^{u_i(s)}sds$. $\sigma=(\sigma_1,..,\sigma_n)\in \Pi$.
Let $$\Pi_1:=\{\beta=(\beta_1,..,\beta_{n-1});\quad (\ref{jan21e1} ) \mbox{ has a solution }\, \}. $$
Note that by Lemma \ref{domain}, $\Pi_1=\mathbb R^{n-1}$ if $a_{ii}>0$ for all $i\in I$.
The mapping from $\Pi_1$ to $\Pi$
is surjective. Here we claim that it is locally one to one. Indeed,
let $\mathbf{M}$ be the following matrix:
$$\mathbf{M}=\left(\begin{array}{ccc}
\partial_{\beta_1}\sigma_1 & \ldots & \partial_{\beta_{n-1}}\sigma_1
\\
\vdots &\vdots & \vdots \\
\partial_{\beta_1}\sigma_{n-1} & \ldots
&\partial_{\beta_{n-1}}\sigma_{n-1}
\end{array}
\right )
$$

We claim that $\mathbf{M}$ is nonsingular
for $\beta\in \Pi_1$ and $\sigma\in \Pi$. We prove this claim by contradiction. Suppose there
exist a non-zero vector $\mathbf{C}=(c_1,..,c_{n-1})^T$ such that
$\mathbf{MC=0}$. Then by setting
$\mathbf{\beta}=c_1\beta_1+...+c_{n-1}\beta_{n-1}$ we have
\begin{equation}\label{jan21e2}
\partial_{\beta}\sigma_1=\partial_{\beta}\sigma_2=...=\partial_{\beta}\sigma_{n-1}=0.
\end{equation}
On the other hand, $\Pi$ is defined by $\Lambda_I=0$, which reads
$$\sum_{i,j\in I}a_{ij}\sigma_i\sigma_j=4\sum_{i\in I}\sigma_i. $$
By differentiating both
sides with respect to $\beta$ we have
$$\sum_i(\sum_ja_{ij}\sigma_j-2)\partial_{\beta}\sigma_i=0. $$
Since $\sum_ja_{ij}\sigma_j>2$, (\ref{jan21e2}) implies
$\partial_{\beta}\sigma_n=0$. Set $\phi_i=\partial_{\beta}u_i$
($i\in I$), then $\phi=(\phi_1,..,\phi_n)$ satisfies the linearized
equation (\ref{jan20e1}) and $\phi_n(0)=0$. From
$\partial_{\beta}\sigma_i=0$ ($i\in I$) we have
$$\int_0^{\infty}e^{u_i}\phi_i(s)sds=0, \quad i\in I$$
which implies from (\ref{jan20e1}) that $\phi$ is bounded at
infinity. By Lemma \ref{jan20lem} $\phi_i=ru_i'+2$, then we see
immediately that this is not possible as $\phi_n(0)=0$. Therefore we
have proved that $\mathbf{M}$ is nonsingular for all
$\beta=(\beta_1,...,\beta_{n-1})\in \Pi_1$.
\medskip

We further assert that there is one-to-one correspondence between $\Pi_1$ and $\Pi$. This is proved in two steps as follows.
\medskip

\noindent Case 1: $a_{ii}>0$, $i\in I$.

In this case, $\Pi_1=\mathbb R^{n-1}$. The mapping from $\Pi_1$ to $\Pi$ is proper and locally one-to-one. Since
both $\mathbb R^{n-1}$ and $\Pi$ are simply connected, there is a one to one correspondence between them.
 Let $u=(u_1,..,u_n)$
and $v=(v_1,..,v_n)$ be two radial solutions such that $u_n(0)=v_n(0)=0$, $\int_{\mathbb R^2}e^{u_i}
=\int_{\mathbb R^2}e^{v_i}$ ($i\in I$). Then $u_i(0)=v_i(0)$ ($i=1,..,n-1 $). Consequently $u_i\equiv v_i$
($i\in I$). Theorem
\ref{thm6} is proved for this case.

\medskip
\noindent Case 2: There exists $i_0\in I$ such that $a_{i_0,i_0}=0$.

We prove this case by a contradiction. Suppose $\beta^k=(\beta^k_1,..,\beta^k_{n-1})\in \Pi_1$ for $k=1,2$ and
$\beta^1\neq \beta^2$, let $u^k$ be the solution corresponding to $\beta^k$ such that $\int_{\mathbb R^2}e^{u^1_i}
=\int_{\mathbb R^2}e^{u^2_i}=\sigma_i$ ($i\in I$).

Just like the case for two equations, we consider the following system
\begin{equation}\label{feb23e1}
\left\{\begin{array}{ll}
u_i''(r)+\frac{u_i'(r)}r+\sum_j(a_{ij}+\epsilon \delta_{ij})e^{u_j}=0\quad 0<r<\infty,\quad i\in I\\
\\
\int_0^{\infty}e^{u_i(r)}rdr<\infty, \quad i\in I. \\   \\
u_1(0)=\beta_1, \,\, ... \,\, ,u_{n-1}(0)=\beta_{n-1},\quad u_n(0)=0.
\end{array}
\right.
\end{equation}

Let $u^{k,\epsilon}$ be the solution to (\ref{feb23e1}) that corresponds to the initial condition $\beta^k$ ($k=1,2$). Let
$\sigma^{k,\epsilon}=(\sigma^{k,\epsilon}_1,..,\sigma^{k,\epsilon}_n)$ be defined as
$\sigma^{k,\epsilon}_i=\int_0^{\infty}re^{u^{k,\epsilon}_i(r)}dr$ ($i=1,..,n$).
By the same argument as in the case of two equations, we have $\sigma^{k,\epsilon}=(\sigma_1,..,\sigma_n)+\circ(1)$ ($k=1,2$) and
$$\frac{\partial \sigma^{k,\epsilon}_i}{\partial \beta_j}=\frac{\partial \sigma_i}{\partial \beta_j}+\circ(1),
\quad i=1,..,n,\quad j=1,..,n-1,\quad k=1,2. $$
Consequently the matrix
$$\left(\begin{array}{ccc}
\partial_{\beta_1}\sigma^{k,\epsilon}_1 & \ldots & \partial_{\beta_{n-1}}\sigma^{k,\epsilon}_1
\\
\vdots &\vdots & \vdots \\
\partial_{\beta_1}\sigma^{k,\epsilon}_{n-1} & \ldots
&\partial_{\beta_{n-1}}\sigma^{k,\epsilon}_{n-1}
\end{array}
\right )
$$
is non-singular at $\beta^1$ or $\beta^2$ for $\epsilon$ small.
On the other hand, $\sigma^{1,\epsilon}$ and $\sigma^{2,\epsilon}$ both satisfy
\begin{equation}\label{feb25e5}
\left\{\begin{array}{ll}
\Lambda^{\epsilon}_I(\sigma^{k,\epsilon})=4\sum_{\in I}\sigma^{k,\epsilon}_i-\sum_{i,j\in I}(a_{ij}+\epsilon \delta_{ij})
\sigma^{k,\epsilon}_i\sigma^{k,\epsilon}_j=0 \\
\\
\Lambda^{\epsilon}_J>0,\quad 0\varsubsetneqq  J\varsubsetneqq  I.
\end{array}
\right.
\end{equation}
We use $\Pi^{\epsilon}$ to represent the hyper-surface described as above.
For
$\sigma^{2,\epsilon}=(\sigma^{2,\epsilon}_1,..,\sigma^{2,\epsilon}_n)\in \Pi^{\epsilon}$, we can find
$\beta^{1,\epsilon}=(\beta^{1,\epsilon}_1,..,\beta^{1,\epsilon}_{n-1})$ such that
$$\beta^{1,\epsilon}_j=\beta^1_j+\circ(1),\quad j=1,2,..,n-1$$ and a solution $\bar u^{1,\epsilon}$ of (\ref{feb23e1})
with the initial condition $(\beta^{1,\epsilon}_1,..,\beta^{1,\epsilon}_{n-1},0)$
such that
$$\int_0^{\infty}re^{\bar u^{1,\epsilon}_j}dr=\sigma^{2,\epsilon}_j,\quad j=1,2,..,n-1. $$
After using $\Lambda_I^{\epsilon}=0$ in (\ref{feb25e5}) we have
$$\int_0^{\infty}re^{\bar u^{1,\epsilon}_n}dr=\sigma^{2,\epsilon}_n. $$
Then the difference
between $\beta^1$ and $\beta^2$ implies $\beta^{1,\epsilon}\neq \beta^2$ for $\epsilon$ small.
A contradiction to the uniqueness property satisfied by the system (\ref{feb23e1}).
Theorem \ref{thm6} is proved for all the cases. $\Box$

\section{Proof of Theorem \ref{thmlocal}}

First we state a Brezis-Merle type Lemma:
\begin{lem}
\label{dec30lem1} Let $\Omega$ be an open, smooth, bounded subset of $\mathbb R^2$. If
$$\sum_j\int_{\Omega}a_{ij}h_j^ke^{u_j^k}\le 4\pi-\delta,\quad i\in
I=\{1,..,n\}$$
 for some $\delta>0$, then for any $\Omega_1\subset\subset
\Omega$, there exists $C(\delta,\Omega,\Omega_1)>0$ such that
$$u_i^k(x)\le C, \quad x\in \Omega_1\subset \subset \Omega,\quad i\in I $$
\end{lem}

\noindent{\bf Proof of Lemma \ref{dec30lem1}:} Let $f_i^k$ ($i\in
I$) be defined as
$$\left\{\begin{array}{ll}
-\Delta f_i^k(x)=\sum_ja_{ij}h_j^ke^{u_j^k},\quad \Omega,\\
\\
f_i^k(x)=0,\quad \mbox{on}\quad \partial \Omega.
\end{array}
\right.
$$
Then by Theorem 1 of \cite{BM}, we have
$$\int_{\Omega}e^{(1+\delta_1)f_i^k}dx\le C,
$$ where $\delta_1>0$ depends on $\delta$.
For any $\Omega'\subset\subset \Omega$, let $x\in \Omega'$, suppose
$B(x,\delta_2)\subset \Omega$, we have, by the mean value property
\begin{eqnarray*}
u_i^k(x)-f_i^k(x)&=&\frac{1}{|B(x,\delta_2)|}\int_{B(x,\delta_2)}(u_i^k(y)-f_i^k(y))dy\\
&\le &C\int_{B(x,\delta_2)}(u_i^k(y)-f_i^k(y))^+dy\\
&\le &C\int_{\Omega}(e^{u_i^k}+e^{f_i^k})\le C,\quad i\in I.
\end{eqnarray*}
So by writing $u_i^k$ as $u_i^k-f_i^k+f_i^k$ we see that
$e^{u_i^k}\in L^{1+\delta_1}(\Omega')$, $i\in I$.  Let $\bar f_i^k$
be defined as
$$\left\{\begin{array}{ll}
-\Delta \bar f_i^k(x)=\sum_{j\in I}a_{ij}h_j^ke^{u_j^k(x)},\quad
\Omega',\\  \\
\bar f_i^k(x)=0,\quad \mbox{on}\quad \partial \Omega'
\quad i\in I.
\end{array}
\right.
$$

Then standard elliptic estimate gives $|\bar f_i^k|\le C$ in
$\Omega'$ ($i\in I$). Let $\Omega''\subset\subset \Omega'$, then for
$x\in \Omega''$, as before we have
$$u_i^k(x)=u_i^k(x)-\bar f_i^k(x)+\bar f_i^k(x)\le
C\int_{\Omega'}(e^{u_i^k}+e^{\bar f_i^k})+C\le C. $$  Lemma
\ref{dec30lem1} is established. $\Box$

\medskip

Recall that $ \sigma_i=\frac 1{2\pi}\int_{\mathbb R^2}h_ie^{v_i}$
($i\in I$) where $h_i=\lim_{k\to \infty}h_i^k(0)$. Since
$v=(v_1,..,v_n)$ satisfies the Liouville system in $\mathbb R^2$, we
have
\begin{equation}\label{jan16e1}
\sum_{j\in I}a_{ij}\sigma_j>2,\quad i\in I.
\end{equation}

Let $\bar \sigma_i=\lim_{r\to 0}\lim_{k\to \infty}\frac 1{2\pi}\int_{B_r}h_i^ke^{u_i^k}$, then the assumption in
Theorem \ref{thmlocal} implies
\begin{equation}\label{feb6e1}
\bar \sigma_i\ge \sigma_i \quad i\in I.
\end{equation}
So (\ref{jan16e1}) also holds for
$\{\bar \sigma_i\}_{i\in I}$.

\begin{lem} \label{dec30lem3}
\begin{equation}
\label{lin2} \sum_{i,j\in I}a_{ij}\bar \sigma_i\bar \sigma_j=4\sum_{i\in
I}\bar \sigma_i.
\end{equation}
\end{lem}

\noindent{\bf Proof of Lemma \ref{dec30lem3}:}

 In the first step we prove that in a small neighborhood of $0$, say, $B(0,r_0)$,
 $u_i^k|_{\partial B_R}\to -\infty$ for $i\in I$ and any fixed
 $0<R<r_0$.

 Indeed, since (\ref{jan16e1}) holds for $\bar \sigma=(\bar \sigma_1,..,\bar \sigma_n)$,
 we have $\sum_{j\in I}a_{ij}\bar \sigma_j>2+3\epsilon_0$ ($i\in I$) for some $\epsilon_0>0$.
 By the definition of $\bar \sigma_i$, we find $r_0$ small and $r_k\to 0$ such that
 $\int_{B_{r_0}\setminus
 B_{r_k}}e^{u_i^k}\le \epsilon_0$ ($i\in I$). Let
 $$\tilde v_i^k(y)=u_i^k(r_ky)+2\ln r_k,\quad |y|\le r_k^{-1}r_0,\quad i\in I.$$
 Then the equation for $\tilde v_i^k$ is
 $$
 -\Delta \tilde v_i^k=\sum_{j\in I}a_{ij}h_j^k(r_k\cdot )e^{\tilde v_j^k},\quad |y|\le
 r_k^{-1}r_0.
$$
Let
$$\bar v_i^k(r)=\frac 1{2\pi r}\int_{\partial B_r}\tilde v_i^k,\quad 1\le
r\le r_k^{-1}r_0,\quad i\in I. $$ Then
$$(\bar v_i^k)'(r)=\frac 1{2\pi r}\int_{B_r}\Delta \tilde v_i^k=-\frac 1{2\pi
r}\int_{B_r}\sum_ja_{ij}h_j^k(r_k\cdot)e^{\tilde v_j^k}dy.$$

For $r>1$,
$$\int_{B_r}\sum_ja_{ij}h_j^k(r_k\cdot)e^{\tilde v_j}>4\pi+2\epsilon_0,\quad
i\in I.$$ So by the definition of $\tilde v_i^k$,
$$(\bar v_i^k)'(r)\le (-2-\frac{\epsilon_0}{\pi})r^{-1},\quad r>1,\quad i\in I.$$
Consequently
$$\bar v_i^k(r_k^{-1}r_0)\le -(2+\frac{\epsilon_0}{\pi})\ln r_k^{-1}+C\to
-\infty,\quad i\in I.
$$
For any fixed $R\in (0,r_0)$,$u_i^k$ has bounded oscillation on any
$\partial B_R$, then we know $u_i^k\to -\infty$ uniformly on
$\partial B_R$. As an immediate consequence, $u^k$ converges to
$-\infty$ on all compact subsets of $B_1\setminus \{0\}$ because
$u^k$ is bounded above in $B_1\setminus B_R$ and $u^k$ has bounded
oscillation on $\partial B_1$.

\medskip

 The second step is to use the first step to evaluate all the terms
in the Pohozaev Identity. Let $G(x,y)$ be the Green's function with
the Dirichlet condition. By the Green's representation formula we
have:
$$u_i^k(x)=\int_{B_1}G(x,y)\sum_ja_{ij}h_j^ke^{u_j^k}-\int_{\partial
B_1}\frac{\partial G(x,y)}{\partial \nu}u_i^k(y)dS_y,\quad i\in I.$$

The Pohozaev identity for the system (\ref{assume1}) defined on
$\Omega$ is of the following form (see the appendix for the proof):
\begin{eqnarray*}
&&\sum_{i\in I}\bigg (\int_{\Omega}(x\cdot \nabla
h_i^k)e^{u_i^k}+2h_i^ke^{u_i^k}\bigg ) \\
&=&\int_{\partial \Omega}\bigg (\sum_i(x\cdot
\nu)h_i^ke^{u_i^k}+\sum_{i,j}a^{ij}\partial_{\nu}u_j^k(x\cdot \nabla
u_i^k)-\frac 12a^{ij}(x\cdot \nu)(\nabla u_i^k\cdot \nabla
u_j^k)\bigg ).
\end{eqnarray*}

Let $\Omega=B_R$ ($R\in (0,1)$) in the Pohozaev Identity, using the
fact that $u_i^k\to -\infty$ in $C^2_{loc}(B_1\setminus \{0\})$ we
observe that
$$\int_{\partial B_R}\sum_i(x\cdot \nu)h_i^ke^{u_i^k}\to 0 \quad
\mbox{and}\quad \int_{B_R}(x\cdot \nabla h_i^k)e^{u_i^k}\to 0, \quad
i\in I.$$ Also we have
$$\frac 1{2\pi}\int_{B_R}\sum_i2h_i^ke^{u_i^k}\to 2\sum_i\bar \sigma_i. $$
For $|x|=R$,
$$\nabla
u_i^k(x)=\int_{B_1}\nabla_xG(x,y)\sum_ja_{ij}h_j^ke^{u_j^k}-
\int_{\partial B_1}\nabla_x(\frac{\partial G(x,y)}{\partial
\nu})u_i^k(y),\quad i\in I$$

The second term of the above is the gradient of a harmonic function
that has bounded oscillation on $\partial B_1$. Let $k\to \infty$,

\begin{equation}\label{9jan4e1}
\partial_ru_i^k(x)\to \frac{\sum_ja_{ij}\bar \sigma_j}{R}+O(1), \quad
\partial_{\theta}u_i^k(x)\to O(1),\,\, i\in I,\,\, |x|=R.
\end{equation}
Using (\ref{9jan4e1}) in the Pohozaev Identity, we have
$$\sum_{ij}a_{ij}\bar \sigma_i\bar \sigma_j=4\sum_i\bar \sigma_i+O(R).
$$
Lemma \ref{dec30lem3} is established by letting $R\to 0$. $\Box$

\medskip

Now we claim
\begin{equation}\label{jan16e2}
\bar \sigma_i=\sigma_i,\quad i\in I.
\end{equation}

To see this, let $s_i=\bar \sigma_i-\sigma_i$. We know from (\ref{feb6e1}) that
$s_i\ge 0$ ($i\in I$). Since for $\{\sigma_i\}_{i\in I}$ we also
have
$$\sum_{ij}a_{ij}\sigma_i\sigma_j=4\sum_i\sigma_i$$
we obtain the following equation for $s_i$ from Lemma \ref{dec30lem3} and the above:
$$\sum_j(\sum_ia_{ij}\bar \sigma_i)s_j+\sum_i(\sum_ja_{ij}\sigma_j)s_i=4\sum_{i}s_i. $$
Since both $\sum_ia_{ij}\bar \sigma_i$ and $\sum_ja_{ij}\sigma_j$
are greater than $2$,
 it is easy to see from the above that $s_i=0$ ($i\in I$).
 (\ref{jan16e2}) is proved.

\medskip

Let $\epsilon_k=e^{-\frac{u_1^k(x_1^k)}2}$, $\bar
h_i^k(y)=h_i^k(\epsilon_ky+x_1^k)$ ($i\in I$). Here we recall that $u_1^k(x_1^k)=\max_{B_1}u_i^k$ ($i\in I$). Then we have
$$
-\Delta v_i^k=\sum_ja_{ij}\bar h_j^ke^{v_j^k}, \quad
\Omega_k, \quad i\in I
$$
where $\Omega_k:=\{y;\,\, \epsilon_ky+x_1^k\in B_1\,\, \}$. Let
\begin{equation}\label{may16e3}
\sigma_i^k=\frac 1{2\pi}\int_{B_1}h_i^ke^{u_i^k},\mbox{ and }
m_i^k=\sum_ja_{ij}\sigma_j^k \quad i\in I.
\end{equation}
 We have $\sigma_i^k\to \sigma_i$ and $m_i^k\to m_i>2$ ($i\in I$).

\begin{prop}\label{1028p1} Given $\delta>0$, there exists
$R(\delta,A,c_0,c_1,\sigma)>1$ such that for all large $k$
\begin{equation}\label{1023e1}
(-m_i^k-\delta)\ln |y|\le v_i^k(y) \le (-m_i^k+\delta)\ln
|y|,\quad y\in \Omega_k\setminus B_{2R},\quad i\in I.
\end{equation}
\end{prop}

\noindent{\bf Proof of Proposition \ref{1028p1}:}

By the convergence of $v_i^k$ to $v_i$ in $C^2_{loc}(\mathbb R^2)$
we only need to prove (\ref{1023e1}) for $2R<|y|\le
\epsilon_k^{-1}$ where $R>>1$.  By the Green's representation formula we have, for $x\in B_1$ and $i\in I$
\begin{equation}\label{26e2}
u_i^k(x)=\int_{B_1}G(x,z)(\sum_ja_{ij}h_j^ke^{u_j^k(z)})-\int_{\partial
B_1}\frac{\partial G(x,z)}{\partial \nu}u_i^k(z).
\end{equation}
 Since the major term of the Green's function is
$-\frac{1}{2\pi}\ln |x-z|$ and the oscillation of $u_i^k$ on
$\partial B_1$ is bounded, we have
$$u_i^k(x)-u_i^k(x_i^k)=\frac 1{2\pi}\int_{B_1}\ln
\frac{|x_i-z|}{|x-z|}(\sum_ja_{ij}h_j^ke^{u_j^k(z)})dz+O(1).
$$
where $u_i^k(x_i^k)=\max_{B_1}u_i^k$. Since our assumption is that $u^k$ converges to $v=(v_1,..,v_n)$
after scaling. The radial symmetry of $v_i$ implies
$$|u_i^k(x_i^k)-u_j^k(x_j^k)|\le C, \quad e^{-\frac 12u_1^k(x_1^k)}|x_i^k-x_j^k|\to 0,\quad i,j\in I. $$
With this observation and the definition of $v_k$ (\ref{26e2}) can be rewritten as
\begin{equation}
\label{eq3} v_i^k(y)=\frac{1}{2\pi}\int_{\Omega_k}\ln
\frac{|z|}{|y-z|}(\sum_ja_{ij}\bar h_j^ke^{v_j^k(z)})dz+O(1),\quad i\in I.
\end{equation}
The proof of (\ref{1023e1}) can be put into two steps. First we
show: For $N>1$, there exists $R>>1$ such that for $|y|>2R$ and
all large $k$,
\begin{equation}\label{eq1}
v_i^k(y)\le -2\ln |y|-N, \quad |y|>2R,\quad i\in I.
\end{equation}

To this end, we use the argument in
Lemma \ref{dec30lem1}. Since $\sigma_i^k\to \sigma_i$, for
$\epsilon>0$ small to be determined, we choose $R>>1$ such that
$$\int_{\Omega_k\setminus B_R}e^{v_i^k}\le \epsilon,\quad i\in I
$$
Fix $r>2R$ and set
$$\bar v_i(z)=v_i^k(rz)+2\ln r+2N,\quad \frac 12<|z|<2,\quad i\in I. $$
By letting $\bar h_i(z)=\bar h_i^k(rz)$ we have
$$
-\Delta \bar v_i(z)=\sum_ja_{ij}\bar h_j(z)e^{-2N}e^{\bar v_j(z)},\quad \frac 12<|z|<2,\, i\in I.
$$
Note that for simplicity we omit $k$ in $\bar v_i(z)$ and $\bar h_i$.
It is readily verified that
$$\int_{\frac 12<|z|<2}e^{\bar v_i(z)}dz\le
e^{2N}\int_{\Omega_k\setminus B_R}e^{v_i^k(y)}dy \quad i\in I.
$$ Now we choose $\epsilon$ to be small enough
so that
$$e^{2N}\int_{\Omega_k\setminus B_R}\sum_ja_{ij}\bar h_j^ke^{v_j^k}\le 3\pi,\quad i\in I. $$
The inequality above implies
\begin{equation}\label{jan26e1}
\int_{B_2\setminus B_{\frac 12}}e^{\bar v_i^k}\le C, \quad i\in I
\end{equation}
where $C$ is independent of $N$.
Using (\ref{jan26e1}) and the argument in Lemma \ref{dec30lem1} we have
\begin{equation}\label{eq2}
\bar v_i(z)\le c_0,\quad |z|=1, \quad i\in I
\end{equation}
 where $c_0$ is a
universal constant. (\ref{eq1}) follows immediately from
(\ref{eq2}).

\medskip

In the second step we use (\ref{eq1}) and (\ref{eq3}) to prove (\ref{1023e1}).
First since $|z|\sim |y-z|$ for $|z|>2|y|$, we have
$$v_i^k(y)=\frac{1}{2\pi}\int_{B_{2|y|}}\ln \frac{|z|}{|y-z|}(\sum_{j}a_{ij}\bar
h_j^ke^{v_j^k(z)})dz+O(1).$$ Next we show
that
\begin{equation}\label{jan16e3}
\frac{1}{2\pi}\int_{B_{2|y|}}|\ln |z||(\sum_ja_{ij}\bar h_j^ke^{v_j^k(z)})dz\le \frac{\delta}{10}\ln |y|,
\quad |y|>R_1
\end{equation}
 where $R_1$ will be chosen large in terms of
$\delta$. Indeed, we can choose $R_1$ so large that
\begin{equation}\label{26e3}
\frac 1{2\pi}\int_{B_{2|y|}\setminus B_{R_1}}\sum_ja_{ij}\bar h_j^ke^{v_j^k(z)}dz<\delta/10.
\end{equation}
 Then
the integral in (\ref{jan16e3}) can be divided into two parts, one part is the integration over $B_{R_1}$, the other
part is the integration on $B_{2|y|}\setminus B_{R_1}$.
Since $e^{v_i}$ decays faster than $|y|^{-2-\delta_1}$ for
some $\delta_1>0$, we use the convergence of $v_i^k$ to $v_i$ to obtain that the integration over
$B_{R_1}$ is $O(1)$. For the other term it is easy to see from (\ref{26e3}) that the integration
over $B_{2|y|}\setminus B_{R_1}$ is less than
$\frac{\delta}5 \ln |y|$. The last term to deal with
is
$$-\frac 1{2\pi}\int_{B_{2|y|}}\ln |y-z|(\sum_ja_{ij}\bar h_j^ke^{v_j^k(z)})dz. $$ For this we divide $B_{2|y|}$ into two
sub-regions:
$$\Omega_1=\{z\in \Omega_k;\,\, |z|<|y|/2\},\quad
\Omega_2:=B_{2|y|}\cap \Omega_k\setminus \Omega_1.$$

Since $|y-z|\sim |y|$ for $z\in \Omega_1$ and
$$|\frac 1{2\pi}\int_{\Omega_1}\sum_ja_{ij}\bar h_je^{v_j^k}-m_1^k|\le \frac{\delta}{20}$$
for $|y|$ large. We obtain immediately that

$$|\frac 1{2\pi}\int_{\Omega_1}\ln |y-z|(\sum_{j}a_{ij}\bar h_i^ke^{v_i^k(z)})dz-m_i^k\ln |y||\le \frac{\delta}{10}\ln |y|,
\quad |y|>R_1.$$
To estimate the last term:
$-\frac 1{2\pi}\int_{\Omega_2}\ln |y-z|(\sum_ja_{ij}\bar h_j^ke^{v_j^k(z)})dz$, we use polar coordinates and (\ref{eq1})
to obtain
$$
|\int_{\Omega_2}\ln |y-z|(\sum_ja_{ij}\bar h_j^ke^{v_j^k(z)})dz|\le Ce^{-N}\ln |y|
$$
for a universal constant $C$.
Choose $N$ large enough we see this
term is less than $\frac{\delta}{10} \ln |y|$. Proposition
\ref{1028p1} is established. $\Box$

\medskip

Since $m_i^k\to m_i>2$, $e^{v_i(y)}\sim O(|y|^{-2-\delta_2})$ for
some $\delta_2>0$. Using this in the proof of Proposition
\ref{1028p1} again we see that
\begin{equation}\label{fortype2}
|v_i^k(y)-m_i^k\ln (1+|y|)|\le C(A,c_0,c_1,\sigma),\quad y\in
\Omega_k,\quad i\in I.
\end{equation}

\begin{prop}\label{pi}
\begin{equation}\label{1030e1}
\sum_{ij}a_{ij}\sigma_i^k\sigma_j^k=4\sum_i\sigma_i^k
+O(\epsilon_k^c)
\end{equation}
where $c>0$ is a small number.
\end{prop}

\begin{rem} Proposition \ref{pi} is equivalent to the second statement of Theorem \ref{thmlocal}.
\end{rem}

\noindent{\bf Proof of Proposition \ref{pi}:} Let $m>2$ be less than $m_i^k$
($i\in I$) and $L_k=\epsilon_k^{-c}$ for $c>0$ small. We estimate
each term of the Pohozaev Identity on $E_k:=B(0,L_k)$:
\begin{eqnarray*}
&&\sum_i\bigg (\int_{E_k}(y\cdot \nabla \bar h_i^k)e^{v_i^k}+2\bar h_i^ke^{v_i^k}\bigg )\\
&=&L_k\int_{\partial E_k}\bigg (\sum_i\bar h_i^ke^{v_i^k}+\sum_{ij}(a^{ij}\partial_{\nu}v_i^k\partial_{\nu}
v_j^k-\frac 12a^{ij}(\nabla v_i^k\cdot \nabla v_j^k))\bigg ).
\end{eqnarray*}
By the decay rate of $v_j^k$ ($j=1,2$), we have
$$
\int_{E_k}(y\cdot \nabla \bar h_i^ke^{v_i^k})=\epsilon_k \int_{E_k}(y\cdot \nabla
h_i^k(\epsilon_ky+x_1^k)e^{v_i^k}
=O(\epsilon_k),\quad i\in I.
$$

$$\int_{E_k}2\bar h_i^ke^{v_i^k}=4\pi \sigma_i^k+O(L_k^{-m+2}), \quad i\in I.
$$
Similarly
$$\int_{\partial
E_k}L_k\bar h_i^ke^{v_i^k}=O(L_k^{-m+2}),\quad i\in I.
$$ Now we estimate $\nabla v_i^k$ ($i\in I$).
By the Green's representation formula:
\begin{eqnarray}
\label{jan3e2} \nabla v_i^k(y)
&=&\int_{\Omega_k}\nabla_yG(y,\eta)(\sum_ja_{ij}\bar
h_j^ke^{v_j^k(\eta)})d\eta \\
&& -\int_{\partial \Omega_k}\nabla_y(\frac{\partial
G(y,\eta)}{\partial \nu})v_i^k(\eta)dS_{\eta} \quad i\in I.\nonumber
\end{eqnarray}
The last term above is the gradient of a harmonic function. We
know that if $f$ is a harmonic function on $B_R$, then $|\nabla
f(0)|\le C\cdot osc(f)/R.$ By this reason we know that, since
$v_i^k$ has bounded oscillation on $\partial \Omega_k$ and
$|y|=L_k<<\epsilon_k^{-1}$, the last term of (\ref{jan3e2}) is
$O(\epsilon_k)$.

To estimate the first term of (\ref{jan3e2}), we use
$$G(y,\eta)=-\frac 1{2\pi}\ln |y-\eta |+H_k(y,\eta). $$
For $|y|=L_k$, $H_k(y,\eta)$, as a function of $\eta$, is a harmonic function of the order
$O(\ln \epsilon_k^{-1})$ on $\partial \Omega_k$. So for $\eta\in E_k$, using $H_k(y,\eta)=H_k(\eta,y)$ and standard
gradient estimate for
harmonic functions, we have
$$|\nabla_yH_k(y,\eta)|=|\nabla_{\eta}H_k(y,\eta)|\le C\frac{\max_{\partial \Omega_k}H_k}{\epsilon_k^{-1}}
=O(\epsilon_k^{\delta}).$$
Consequently
$$\int_{E_k}\nabla_yH_k(y,\eta)(\sum_ja_{ij}\bar h_j^ke^{v_j^k})=O(\epsilon_k^{\delta}) $$
for $\delta\in (0,1)$.
We are left with the estimate of the term
$$-\frac 1{2\pi}\int_{E_k}\nabla_y(\ln |y-\eta |)\sum_ja_{ij}\bar
h_j^ke^{v_j^k(\eta)}d\eta. $$ For this we use
\begin{eqnarray*}
&&\partial_{y_a}(-\frac 1{2\pi}\ln
|y-\eta|)-\partial_a(-\frac{1}{2\pi}\ln |y|)\\
&=&-\frac 1{2\pi}\frac{-\eta_a|y|^2-y_a|\eta
|^2+2y_a\sum_{t=1}^2y_t\eta_t}{|y-\eta|^2|y|^2},\quad a=1,2
\end{eqnarray*}
and elementary estimate to obtain
\begin{eqnarray*}
&&-\frac 1{2\pi}\int_{E_k}\nabla_y(\ln |y-\eta |-\ln
|y|)(\sum_ja_{ij}\bar
h_j^ke^{v_j^k(\eta)})d\eta\\
&=&O(L_k^{-m+1}\ln L_k).
\end{eqnarray*}
 Consequently
\begin{eqnarray*}
\partial_av_i^k(y)&=&\int_{\Omega_k}\partial_a(-\frac 1{2\pi}\ln
|y|)(\sum_ja_{ij}\bar h_j^ke^{v_j^k})d\eta\\
&=&-m_i^k\frac{y_a}{|y|^2}+O(L_k^{-m+1}\ln L_k),\quad i\in
I,\,\, a=1,2.
\end{eqnarray*}

Using this in the computation of the Pohozaev Identity we obtain
(\ref{1030e1}). Proposition \ref{pi} is established. $\Box$

\bigskip

Now we are in the position to prove (\ref{uniformes}).
 One can find $\{\sigma_{i,k}\}_{i\in I}$ that
satisfies $\Lambda_I(\sigma_{\cdot,k})=0$, which is
$$\sum_{i,j}a_{ij}\sigma_{ik}\sigma_{jk}=4\sum_i\sigma_{ik}$$
so that
\begin{equation}\label{march2e1}
\sigma_{i,k}=\sigma_i^k, \, i=1,..,n-1,\quad
\sigma_{n,k}-\sigma_n^k=O(\epsilon_k^{\delta})
\end{equation}
 for some $\delta>0$. For $\{\sigma_{i,k}\}_{i\in I}$ we let $\bar V^k=(\bar V^k_1,..,\bar V^k_n)$
be the unique global solution so that $\{\bar V_i^k\}_{i\in I}$ are radial with respect to the origin,
$$\frac 1{2\pi}\int_{\mathbb R^2}h_i^k(0)e^{\bar V_i^k}=\sigma_{i,k}, \quad i\in I, \quad \bar V_1^k(0)=0.$$ Note
that the uniqueness is proved in Theorem \ref{thm6}. Using
$\sigma_{i,k}\to \sigma_i$ ($i\in I$) as $k\to \infty$, we assert that
$\bar V_i^k\to v_i$ ($i\in I$) in $C^2_{loc}(\mathbb R^2)$ because $v=(v_1,..,v_n)$
is the only radial solution that satisfies $\frac 1{2\pi}\int_{\mathbb R^2}h_ie^{v_i}=\sigma_i$ and
$v_1(0)=0$. On the other hand, by standard potential analysis
$$|\bar V_i(y)+\bar m_{i,k}\ln |y||\le C(A,\sigma),\quad
|y|>2$$ where $\bar m_{i,k}=\sum_ja_{ij}\sigma_{j,k}$. (\ref{march2e1}) implies
$|\bar m_{i,k}-m_i^k|=O(\epsilon_k^{\delta})$. Thus by (\ref{fortype2}) we have
$$|v_i^k(y)-\bar V_i^k(y)|\le C(A,c_0,c_1,\sigma), \quad y\in \Omega_k. $$
Let $V_i^k$ be defined by
$$V_i^k(\epsilon_ky)+2\log \epsilon_k=\bar V_i^k(y),$$
then the second statement of Theorem \ref{thmlocal} is established.
$\Box$

\section{Proof of Theorem \ref{location} and (\ref{loca})}

In this section we prove Theorem \ref{location} and (\ref{loca}).
Let
$$\tilde u_i^k=u_i^k-\psi_i^k,\quad \tilde
h_i^k=h_i^ke^{\psi_i^k}, \quad i\in I.$$ Since $\psi_i^k(0)=0$ we
have
$$\frac{\nabla \tilde h_i^k(0)}{\tilde h_i^k(0)}=\nabla
\psi_i^k(0)+\frac{\nabla h_i^k(0)}{h_i^k(0)}. $$ Let $|\xi |=1$ be a
unit vector, then a Pohozaev identity for $\tilde u^k=(\tilde u_1^k,..,\tilde u^k_n)$ is of the
form (see the appendix for the proof)
\begin{eqnarray*}
&&\int_{B_R}(\sum_i\partial_{\xi}\tilde h_i^ke^{\tilde u_i^k})\nonumber\\
&=&\int_{\partial B_R}\bigg ( \sum_ie^{\tilde u_i^k}\tilde
h_i^k(\xi\cdot \nu)+\sum_{ij}a^{ij}(\partial_{\nu}\tilde
u_i^k\partial_{\xi}\tilde u_j^k-\frac 12(\xi\cdot \nu)(\nabla \tilde
u_i^k\cdot \nabla \tilde u_j^k))\bigg )
\end{eqnarray*}
By choosing $0<R<1$, it is easy to see from the decay rate of
$\tilde u_i^k$ that
$$\int_{\partial B_R}\sum_i(e^{\tilde
u_i^k}\tilde h_i^k (\xi\cdot \nu))\to 0.
$$ Also, since $\tilde h_i^ke^{\tilde u_i^k}\to 2\pi\sigma_i\delta_0$ in
distributional sense, the left hand side of the Pohozaev identity
tends to
$$2\pi\sum_i\frac{\partial_{\xi}\tilde h_i(0)}{\tilde
h_i(0)}\sigma_i. $$ To consider the limit of $\nabla \tilde
u_i^k(x)$ for $|x|=R$, we use the Green's representation formula:
$$\tilde u_i^k(x)=\int_{B_1}G(x,\eta)(\sum_ja_{ij}\tilde h_j^k(\eta)e^{\tilde u_j^k(\eta)})+constant. $$ By taking the
derivative on $x$ and letting $k\to \infty$, we have
$$\nabla \tilde u_i^k(x)\to 2\pi
m_i\nabla_1G(x,0)=m_i\frac{x}{|x|}. \quad i\in I.$$  Using this in
the computation of the Pohozaev identity we see the limit of the
right hand side is $0$. Therefore we have obtained:
$$\sum_i\frac{\partial_{\xi}\tilde h_i(0)}{\tilde
h_i(0)}\sigma_i=0.$$ Since $\xi$ is arbitrary, Theorem
\ref{location} is established. $\Box $

\medskip

\noindent{\bf Proof of (\ref{loca}):} Since $\int_Mh_i^ke^{u_i^k}=1$
($i\in I$) the equation for $\{u^k\}$ is (see (\ref{1218e1}))
\begin{equation}\label{jan18e1}
\Delta_g u_i^k+\sum_{j=1}^n\rho_ja_{ij}(h_j^ke^{u_j^k}-1)=0,\quad M.
\end{equation}

 Recall that $\{p_1,..,p_m\}$
are disjoint blowup points for $\{u^k\}$. Let
\begin{equation}\label{jan19e3}
\sigma_{it}=\lim_{r\to 0}\lim_{k\to \infty}\frac
1{2\pi}\int_{B(p_t,r)}h_i^ke^{u_i^k}dV_g.
\end{equation}
 Our assumption is that
around each $p_t$, $\{u^k\}$ converges to a Liouville system of $n$
equations after scaling. Let $\delta>0$ be small enough so that
$B(p_t,\delta)$ ($t=1,..,m$) are disjoint. For each $t$, let $M_t^k$
be the maximum of $\{u_i^k\}_{i\in I}$ in $B(p_t,\delta)$. In the
isothermal coordinates around $p_t$, $g=e^{\phi}\delta_0$ and $\Delta_g=e^{-\phi}\Delta$
where $\delta_0$ is the Euclidean metric. We also have $\phi(0)=|\nabla \phi(0)|=0$.
With these properties (\ref{jan18e1}) in $B(p_t,\delta)$ becomes
$$
\Delta
u_i^k+\sum_{j=1}^n\rho_ja_{ij}e^{\phi}h_j(e^{u_j^k}-1)=0,\quad
B_{\delta},\quad i\in I.
$$
Let $f_i$ satisfy
$$\Delta f_i=\sum_j\rho_ja_{ij}e^{\phi}h_j\quad B_{\delta},\quad i\in I$$
and $f_i=0$ on $\partial B_{\delta}$, then the equation for $u_i^k$
can further be written as
\begin{equation}\label{jan19e1}
\Delta
(u_i^k+f_i)+\sum_j\rho_ja_{ij}e^{\phi-f_i}h_je^{u_j^k+f_i}=0,\quad
B_{\delta},\quad i\in I.
\end{equation}
Let
$$\sigma_{it}^k=\frac
1{2\pi}\int_{B_{\delta}}h_i^ke^{\phi}e^{u_i^k},\quad
m_{it}^k=\sum_ja_{ij}\sigma_{jt}^k.
$$
By Theorem \ref{thmlocal} and $\phi(0)=0$ the limit of
$\sigma_{it}^k$ is $\sigma_{it}$ (defined in (\ref{jan19e3})). Let
$m_{it}>2$ be the limit of $m_{it}^k$, then from Theorem
\ref{thmlocal} we have, for $x\in
\partial B(p_t,\delta)$
\begin{equation}\label{jan19e2}
u_i^k(x)=-\frac{m_{it}^k-2}2M_t^k+O(1),\quad x\in \partial
B(p_t,\delta)\quad i\in I,\quad t=1,..,m.
\end{equation}
From the Green's representation of $u_i^k$ it is easy to see that
the difference between $u_i^k(x)$ and $u_i^k(y)$ for $x,y$ away from
the blowup set is uniformly bounded. Therefore for fixed $t_1$ and $t_2$, using
$M_t^k\to \infty$ we obtain from (\ref{jan19e2}) that
\begin{equation}\label{jan19e4}
\frac{m_{it_1}-2}{m_{it_2}-2}=\lambda_{t_1t_2},\quad i\in I.
\end{equation}
We claim that $\lambda_{t_1t_2}=1$. Indeed, $\{\sigma_{it}\}_{i\in
I}$ satisfies
$$\sum_{ij}a_{ij}\sigma_{it}\sigma_{jt}=4\sum_i\sigma_{it}$$
which can be written as
$$\sum_{ij}a^{ij}m_{it}m_{jt}=4\sum_{ij}a^{ij}m_{jt}. $$
The above is equivalent to
$$\sum_{ij}a^{ij}(m_{it}-2)(m_{jt}-2)=4\sum_{ij}a^{ij}. $$
Replacing $m_{it}$ by $m_{it_1}$ and $m_{it_2}$ respectively in the above, we
have
$$(1-\lambda_{t_1t_2}^2)\sum_{ij}a^{ij}=0. $$
Recall that $A$ is assumed to be positive definite. So $\sum_{ij}a^{ij}>0$, we have
$\lambda_{t_1t_2}=1$ ($t_1,t_2=1,..,m$).

We can further claim that
\begin{equation}\label{jan19e5}
\sigma_{it}=\frac 1{2\pi m},\quad i\in I, \quad t=1,..,m.
\end{equation}
because $\int_{M}h_i^ke^{u_i^k}\equiv 1$ ($i\in I$),
$m_{it_1}=m_{it_2}$ ($i\in I$) and
$$\int_{M\setminus
\cup_{t=1}^mB(p_t,\delta)}e^{u_i^k}dV_g\to 0,\quad i\in I. $$

The Green's representation for $u_i^k$ is
\begin{equation}\label{feb26e1}
u_i^k(x)=\bar
u_i^k+\int_MG(x,\eta)\sum_j\rho_ja_{ij}h_je^{u_j^k}dV_g.
\end{equation}
 The last
term of the above tends to
\begin{equation}\label{feb26e2}
\sum_{t=1}^mG(x,p_t)(\sum_j\rho_ja_{ij})/m.
\end{equation}
Recall that
\begin{equation}\label{feb26e3}
G(x,\eta)=-\frac 1{2\pi}\chi \ln d(x,\eta)+G^*(x,\eta).
\end{equation}
For $x\in \partial B(p_s,\delta)$, by choosing the support of $\chi$
possibly smaller, we observe that $G(x,p_t)=G^*(x,p_t)$ for $t\neq
s$. Therefore, let $\phi_k$ be the harmonic function on
$B(p_s,\delta)$ defined by the oscillation of $u_i^k$ on $\partial
B(p_s,\delta)$, using (\ref{feb26e1}), (\ref{feb26e2}) and (\ref{feb26e3}) we have
$$\lim_{k\to
\infty}\nabla_g\phi_k(p_s)=\sum_{t=1}^m\nabla_1G^*(p_s,p_t)(\sum_j\rho_ja_{ij})/m.
$$
Then (\ref{loca}) is a consequence of Theorem \ref{location} and the
above.

\section{Appendix: The Pohozaev identity for the Liouville system}

In this section we derive the Pohozaev identity for the Liouville
system
\begin{equation}
\label{1219e1} -\Delta u_i=\sum_{j=1}^na_{ij}h_je^{u_j},\quad
\Omega\subset \subset \mathbb R^2, \quad i\in I.
\end{equation}

The Pohozave identity for (\ref{1219e1}) is
\begin{eqnarray}
&&\sum_{i\in I}\bigg (\int_{\Omega}(x\cdot \nabla
h_i)e^{u_i}+2h_ie^{u_i}\bigg )\label{pisys2} \\
&=&\int_{\partial \Omega}\bigg (\sum_i(x\cdot
\nu)h_ie^{u_i}+\sum_{i,j}a^{ij}(\partial_{\nu}u_j(x\cdot \nabla
u_i)-\frac 12(x\cdot \nu)(\nabla u_i\cdot \nabla u_j))\bigg
).\nonumber
\end{eqnarray}

\noindent{\bf Proof of (\ref{pisys2}):} We write (\ref{1219e1}) as
\begin{equation}\label{1219e2}
-\sum_ja^{ij}\Delta u_j=h_ie^{u_i}, \quad \Omega,\quad i\in I.
\end{equation}
 By multiplying $x\cdot \nabla u_i$ to the right hand
side of (\ref{1219e2}) and integration by parts, we obtain the
following terms:
$$\int_{\Omega}(x\cdot
\nu)h_ie^{u_i}-2\int_{\Omega}h_ie^{u_i}-\int_{\Omega}(x\cdot \nabla
h_i)e^{u_i}.$$

Multiply $x\cdot \nabla u_i$ to the left hand side of
(\ref{1219e2}) and use integration by parts, we have, after
taking the summation on $i$
\begin{eqnarray*}
-\sum_{ij}\int_{\partial \Omega}a^{ij}\partial_{\nu}u_jx\cdot \nabla
u_i+\int_{\Omega}\sum_{ij}a^{ij}\nabla u_i\nabla u_j\\
+\sum_{ij}\int_{\Omega}\sum_{a=1}^2\sum_{b=1}^2 a^{ij}x_b\partial_au_j\partial_{ab}u_i.
\end{eqnarray*}
Using the symmetry of $a^{ij}$ and integration by parts again the left hand side is equal to
$$-\sum_{ij}\int_{\partial \Omega}a^{ij}\partial_{\nu}u_jx\cdot \nabla
u_i+\frac 12\sum_{ij}\int_{\partial \Omega}a^{ij}(x\cdot \nu)(\nabla
u_i\cdot \nabla u_j).$$ Then (\ref{pisys2}) follows.

\medskip

A different version of the Pohozaev identity is as follows. Let
$\xi$ be a unit vector, then we have

\begin{eqnarray}\label{pisys3}
&&\sum_i\int_{\Omega}\partial_{\xi}h_ie^{u_i}\\
&=&\int_{\partial \Omega}\sum_ie^{u_i}h_i(\xi\cdot
\nu)+\sum_{i,j}a^{ij}\bigg
(\partial_{\nu}u_i\partial_{\xi}u_j-\frac 12(\xi\cdot \nu)(\nabla
u_i\cdot \nabla u_j)\bigg ).\nonumber
\end{eqnarray}

The third Pohozaev identity is for the linearized system:
$$(r\phi_i'(r))'+\sum_ja_{ij}e^{u_j}r\phi_j(r)=0,\quad 0<r<\infty,\quad i\in I. $$
The Pohozaev identity is:
\begin{equation}\label{pisyslinear}
\sum_i(r^2\phi_i(r)e^{u_i}-2\int_0^rse^{u_i}\phi_ids)=-\sum_{i,j}a^{ij}(r\phi_j'(r))(ru_i'(r)).
\end{equation}
To derive (\ref{pisyslinear}) we just need to write the linear
system as
$$-\sum_ja^{ij}(r\phi_j'(r))'=e^{u_i}\phi_i(r)r,\quad i\in I. $$
Multiply $ru_i'(r)$ to both sides of the above and use integration
by parts, we obtain (\ref{pisyslinear}).


\begin{thebibliography}{99}

\bibitem{aly} J. J. Aly, Thermodynamics of a two-dimensional
self-gravitating system, Phy. Rev. A 49, No. 5, Part A (1994),
3771-3783.

\bibitem{bennet} W. H. Bennet, Magnetically self-focusing streams, Phys. Rev. 45 (1934),
890–897.

\bibitem{bclt} D. Bartolucci, C.C. Chen, C.S. Lin, G. Tarantello, Profile of blow-up solutions
to mean field equations with singular data. Comm. Partial Differential Equations 29 (2004), no. 7-8, 1241--1265.

\bibitem{biler} P. Biler and T. Nadzieja, Existence and nonexistence of solutions of a
model of gravitational interactions of particles I \& II, Colloq.
Math. 66 (1994), 319--334; Colloq. Math. 67 (1994), 297--309.

\bibitem{BM} H. Brezis, F. Merle,Uniform estimates and blow-up behavior for solutions
of $-\Delta u=V(x)e\sp u$ in two dimensions. Comm. Partial
Differential Equations 16 (1991), no. 8-9, 1223--1253.

\bibitem{CGS} L. A. Caffarelli, B. Gidas, J. Spruck,
Asymptotic symmetry and local behavior of semilinear elliptic
equations with critical Sobolev growth. Comm. Pure Appl. Math. 42
(1989), no. 3, 271--297.

\bibitem{chang} S. A. Chang, P. Yang, Conformal deformations of
metrics on $\mathbb S^2$, J. Diff. Geom. 27 (1988), 256-296.

\bibitem{chanillo1} S. Chanillo, M. K-H Kiessling, Rotational symmetry of solutions of some
nonlinear problems in statistical mechanics and in geometry. Comm. Math. Phys. 160 (1994), no. 2, 217--238.

\bibitem{chanillo2}S. Chanillo, M. K-H Kiessling,
Conformally invariant systems of nonlinear PDE of Liouville type.
Geom. Funct. Anal. 5 (1995), no. 6, 924--947.

\bibitem{ChenLin1} C. C. Chen, C. S. Lin, Sharp estimates for solutions of multi-bubbles
in compact Riemann surfaces. Comm. Pure Appl. Math. 55 (2002), no. 6, 728--771.

\bibitem{chenlin2} C. C. Chen, C. S. Lin, Topological degree for a mean field equation on Riemann surfaces.
  Comm. Pure Appl. Math.  56  (2003),  no. 12, 1667--1727.

\bibitem{chenli} W. X. Chen, C.M. Li,Classification of solutions of some nonlinear
elliptic equations. Duke Math. J. 63 (1991), no. 3, 615--622.

\bibitem{chenli2} W. X. Chen, C. M. Li, Qualitative properties of solutions to some nonlinear elliptic
equations in $R\sp 2$. Duke Math. J. 71 (1993), no. 2, 427--439.

\bibitem{childress} S. Childress and J. K. Percus, Nonlinear aspects of Chemotaxis, Math. Biosci. 56 (1981),
217–237.

\bibitem{CSW} M. Chipot, I. Shafrir, G. Wolansky,
On the solutions of Liouville systems. J. Differential Equations
140 (1997), no. 1, 59--105.

\bibitem{CSW1} M. Chipot, I. Shafrir, G. Wolansky,
Erratum: "On the solutions of Liouville systems" [J. Differential
Equations 140 (1997), no. 1, 59--105]; J. Differential Equations
178 (2002), no. 2.

\bibitem{debye} P. Debye and E. Huckel, Zur Theorie der Electrolyte, Phys. Zft 24 (1923), 305–325.

\bibitem{dunne} G. Dunne, Self-dual Chern-Simons Theories, Lecture
Notes in Physics, vol. m36, Berlin: Springer-Verlag, 1995.

\bibitem{jw1} J. Jost, G. Wang, Classification of solutions of a Toda system in ${\mathbb R}\sp 2$.
 Int. Math. Res. Not. 2002, no. 6, 277--290.

\bibitem{jw2} J. Jost, G. Wang, Analytic aspects of the Toda
system. I. A Moser-Trudinger inequality. Comm. Pure Appl. Math. 54
(2001), no. 11, 1289--1319.

\bibitem{jlw} J. Jost, C. S. Lin, G. Wang, Analytic aspects of the Toda
 system. II. Bubbling behavior and existence of solutions.
 Comm. Pure Appl. Math. 59 (2006), no. 4, 526--558.

\bibitem{keller} E. F. Keller and L. A. Segel, Traveling bands of Chemotactic Bacteria: A theoretical analysis,
J. Theor. Biol. 30 (1971), 235–248.

\bibitem{kiessling2}M. K.-H. Kiessling and J. L. Lebowitz, Dissipative stationary Plasmas: Kinetic Modeling Bennet
Pinch, and generalizations, Phys. Plasmas 1 (1994), 1841–1849.

\bibitem{licmp} Y. Y. Li,
Harnack type inequality: the method of moving planes. Comm. Math.
Phys. 200 (1999), no. 2, 421--444.

\bibitem{lishafrir} Y. Y. Li, I. Shafrir, Blow up analysis for solutions
of $ -\Delta u = Ve^u$ in dimension two, Indiana Univ. Math. J. 43
(1994), 1255-1270.

\bibitem{lincla} C. S. Lin, A classification of solutions of a conformally invariant fourth order equation
in $\mathbb R^n$. (English summary) Comment. Math. Helv. 73
(1998), no. 2, 206--231.

\bibitem{linzhang2} C. S. Lin, L. Zhang, Topological degree for
some Liouville systems on Riemann surfaces, in preparation.

\bibitem{malchiodi1} A. Malchiodi, A. Malchiodi,
Morse theory and a scalar field equation on compact surfaces,
preprint.

\bibitem{malchiodi} A. Malchiodi, C. B. Ndiaye, Some existence results for the Toda system on closed surfaces.
 Atti Accad. Naz. Lincei Cl. Sci. Fis. Mat. Natur. Rend. Lincei (9) Mat. Appl. 18 (2007), no. 4, 391--412.

\bibitem{mock} M. S. Mock, Asymptotic behavior of solutions of transport equations for semiconductor devices,
J. Math. Anal. Appl. 49 (1975), 215–225.

\bibitem{rubinstein} I. Rubinstein, Electro diffusion of Ions,
SIAM, Stud. Appl. Math. 11 (1990).

\bibitem{wolansky1} G. Wolansky, On steady distributions
of self-attracting clusters under friction and fluctuations, Arch. Rational Mech. Anal. 119 (1992), 355--391.

\bibitem{wolansky2} G. Wolansky, On the evolution of self-interacting clusters and applications to semi-linear
equations with exponential nonlinearity, J. Anal. Math. 59 (1992), 251--272.

\bibitem{yang} Y. Yang, Solitons in field theory and nonlinear
analysis, Springer-Verlag, 2001.

\bibitem{zhangcmp} L. Zhang, Blowup solutions of some nonlinear elliptic equations
involving exponential nonlinearities. (English summary)
Comm. Math. Phys. 268 (2006), no. 1, 105--133.

\bibitem{zhangccm} L. Zhang, Asymptotic behavior of blowup solutions for elliptic equations
with exponential nonlinearity and singular data, to appear in
Communications on Contemporary mathematics, available at arXiv:
0810.5143

\end{thebibliography}
\end{document}